\documentclass[a4paper,12pt]{article}
\usepackage{amssymb}
\usepackage{graphicx}

\newcommand{\R}{\mathbb{R}}

\newcommand{\N}{\mathbb{N}}

\newcommand{\beq}{\begin{equation} }
\newcommand{\eqq}{\end{equation} }
\newcommand{\cuad}{{\sqcap\kern-.68em\sqcup}}

\newcommand{\norm}[1]{\|#1\|}
\newcommand{\equ}[1]{(\ref{#1})}
\newcommand\rn{\mathbb{R}^N}

\newtheorem{teo}{Theorem}[section]
\newtheorem{prop}{Proposition}[section]
\newtheorem{proposition}{Proposition}[section]
\newtheorem{lema}{Lemma}[section]
\newtheorem{lemma}{Lemma}[section]

\newtheorem{remark}{Remark}[section]
\newcommand{\bremark}{\begin{remark} \em}
\newcommand{\eremark}{\end{remark} }

\newcommand{\tplus}{t^*_+}
\newcommand{\tminus}{t^*_-}
\newcommand{\tl}{t^*_\lambda}
\newcommand{\phiplus}{\varphi^+_1}
\newcommand{\phiminus}{\varphi^-_1}
\newcommand{\F}[1]{F[#1]}
\newcommand{\lambdaplus}{\lambda^+_1}
\newcommand{\lambdaminus}{\lambda^-_1}
\newcommand{ \vae}{\varepsilon}
\newcommand{ \ep}{\varepsilon}

\def\linf{L^\infty(\Omega)}
\def\lp{L^p(\Omega)}

\def\mp{{{\cal M}^+_{\lambda,\Lambda}}}
\def\mm{{{\cal M}^-_{\lambda,\Lambda}}}

\newcommand{\lpl}{\lambda_1^+}
\newcommand{\lmin}{\lambda_1^-}

\newcommand{\fipl}{\varphi_1^+}
\newcommand{\fimin}{\varphi_1^-}

\def\beeq{\begin{equation}}
\def\eeq{\end{equation}}
\newcommand{\begeqaet}{\begin{eqnarray*}}
\newcommand{\eneqaet}{\end{eqnarray*}}
\newcommand{\calc}{{\cal C}}
\newcommand{\cals}{{\cal S}}

\def\co{C(\overline{\Omega})}

\newcommand{\re}[1]{(\ref{#1})}
\newcommand{\hip}{($F^{\ell}_+$)}
\newcommand{\him}{($F^{\ell}_-$)}
\newcommand{\hdp}{($F^{r}_+$)}
\newcommand{\hdm}{($F^{r}_-$)}

\begin{document}
\begin{center}{\bf\Large  Landesman-Lazer type results for second order Hamilton-Jacobi-Bellman equations }\medskip

{Patricio FELMER}, {Alexander QUAAS}, {Boyan SIRAKOV}\end{center}
\setcounter{equation}{0}
\section{ Introduction}

We study the boundary-value  problem
\beeq\label{first}
\left\{\begin{array}{rclcc} F(D^2u,Du,u,x)+\lambda u&=&f(x,u)&\mbox{in}& \Omega,\\
u&=&0&\mbox{on}& \partial\Omega, \end{array}\right.\eeq where the  second
order differential operator $F$ is of Hamilton-Jacobi-Bellman (HJB)
type, that is, $F$ is a supremum of linear elliptic operators, $f$ is sub-linear in $u$ at infinity, and $\Omega\subset\rn$ is a regular bounded domain.

HJB operators have been the object of intensive study during the last thirty years -- for a general review of their theory and  applications we refer to \cite{FS}, \cite{K1}, \cite{Son}, \cite{Ca}. Well-known  examples include the Fucik operator $\Delta u + b u^+ + au^-$ (\cite{Fu}), the Barenblatt operator $\max\{a\Delta u, b\Delta u\}$ (\cite{BER}, \cite{KPV}), and the Pucci operator $\mp(D^2u)$ (\cite{pucci}, \cite{CC}).

To introduce the problem we are interested in, let us first recall some classical results in the case when $F$ is the Laplacian and $\lambda\in (-\infty,\lambda_2)$, (we shall denote with  $\lambda_i$ the $i$-th eigenvalue of the Laplacian). If $f$ is independent of $u$ the solvability of \re{first} is  a consequence of the Fredholm alternative, namely, if $\lambda\not=\lambda_1$, problem \re{first} has a solution for each $f$, while if  $\lambda=\lambda_1$ (resonance) it has solutions if and only if $f$ is orthogonal to $\varphi_1$,
the first eigenfunction of the Laplacian. The existence result in the non-resonant case extends to nonlinearities $f(x,u)$ which grow sub-linearly in $u$ at infinity, thanks to Krasnoselski-Leray-Schauder degree and fixed point theory, see  \cite{Am2}.

A fundamental result, obtained by  Landesman and Lazer \cite{LL}
(see also \cite{Hess}), states that in the resonance case $\lambda=
\lambda_1$ the
problem
$$
\Delta u +\lambda_1 u = f(x,u)\quad\mbox{in}\,\,\Omega, \qquad u=0\quad
\mbox{on}\,\,\partial\Omega,
$$
is solvable provided $f$ is bounded and, setting
\begin{equation}\label{ff}
{f}^\pm(x):=\displaystyle\limsup_{s\to\pm\infty} f(x,s),\quad {f}_\pm(x):=\displaystyle\liminf_{s\to\pm\infty} f(x,s),
\end{equation}
(this notation will be kept from now on) one of the following conditions is satisfied
\begin{equation}\label{ll}
\int_\Omega {f}^-\varphi_1<0<\int_\Omega {f}_+\varphi_1\,,\qquad
\int_\Omega {f}_-\varphi_1>0>\int_\Omega {f}_+\varphi_1.\end{equation}
This result initiated a huge amount of work on solvability of boundary value
problems in which the elliptic operator is at, or more generally close to,
resonance. Various extensions of the results in \cite{LL} for resonant problems
 were obtained in \cite{AAM}, \cite{AM} and \cite{BN}. Further,
Mawhin-Schmidt \cite{MS} -- see also \cite{CMN}, \cite{CdF} --
considered \re{first} with $F=\Delta$ for $\lambda$ close to $ \lambda_1$,
and showed that the first (resp. the second) condition in \re{ll} implies
that for some $\delta>0$ problem \re{first} has at least one solution for
$\lambda \in (\lambda_1-\delta,\lambda_1]$ and at least three solutions for $\lambda \in (\lambda_1,\lambda_1+\delta)$ (resp. at least one solution for $\lambda \in [\lambda_1,\lambda_1+\delta)$ and at least three solutions for $\lambda \in (\lambda_1,\lambda_1+\delta)$). These results rely on degree theory and,
more specifically, on the notion of bifurcation from infinity, studied by Rabinowitz in \cite{R}.

The same results naturally hold if the Laplacian is replaced by any
uniformly elliptic operator in divergence form. Further, they do remain true
if a general linear operator in non-divergence form
\beq\label{linear}
 \displaystyle L=
 a_{ij}(x)\partial^2_{
x_i x_j}+  b_i(x)\partial_{ x_i} +c(x),
\eeq
is considered, but we have to change  $\varphi_1$ in  $\re{ll}$
by the first eigenfunction of the adjoint operator of $L$. This fact is probably known to the experts, though we are not aware of a reference. Its proof -- which will also easily follow from our arguments below -- uses the Donsker-Varadhan (\cite{DV}) characterization of the first eigenvalue of $L$ and the results in  \cite{BNV} which link the positivity of this eigenvalue to the validity of the maximum principle and to the Alexandrov-Bakelman-Pucci inequality (the degree theory argument remains the same as in the divergence case).

The interest in this type of problems has remained high in the PDE community over the years. Recently a large number of works have considered the extensions of the above results to quasilinear equations (for instance, replacing the Laplacian by the $p$-Laplacian), where somewhat different phenomena take place,
see
 \cite{AO},  \cite{delpino}, \cite{DGT1}, \cite{DGT2}.  There has also been a considerable interest in refining the Landesman-Lazer hypotheses \re{ll} and finding general hypotheses on the nonlinearity which permit to determine on which side of the first eigenvalue the bifurcation from infinity takes place, see \cite{AmA}, \cite{AG}, \cite{GR}.

It is our goal here to study the  boundary value problem \re{first} under
Landesman-Lazer conditions on $f$, when $F$ is a Hamilton-Jacobi-Bellman (HJB)
operator, that is, the supremum of linear operators as in \re{linear}:
\beq\label{formf}
F[u]:= F(D^2u, Du, u, x) =
\displaystyle\sup_{\alpha\in {\cal A}} \{\mathrm{tr} (A^\alpha(x)D^2u) + b^\alpha(x).Du
+ c^\alpha(x) u \},
 \eeq
where ${\cal A}$ is an arbitrary index set. The following
hypotheses on $F$  will be kept throughout the paper:
$A^\alpha \in
C(\overline{\Omega})$, $b^\alpha, c^\alpha\in L^\infty(\Omega)$  for
all $\alpha\in {\cal A}$ and, for some constants $0<\lambda\le \Lambda$, we have $\lambda I\! \le\! A^\alpha(x) \le\! \Lambda I$,
for all $
x\in \Omega$ and all $\alpha\in {\cal A}$.
We stress however that all our results are new even for operators with smooth coefficients.

Let us now describe the most distinctive features of HJB operators -- with respect to the operators considered in the previous papers on Landesman-Lazer
type problems -- which make our work and results different.
The HJB operator $F[u]$ defined in \re{formf} is nonlinear, yet positively homogeneous, (that is $F[tu]=tF[u]$ for $t\ge0$), thus  one may
expect it has eigenvalues and eigenfunctions on the cones of positive and negative functions,
but they may be different to each other. This fact was established by Lions
in 1981, in the case of operators with regular coefficients, see \cite{L}.
In that paper he proved $F[u]$ has two real "demi"- or "half"-principal
eigenvalues $\lpl, \lmin\in\mathbb{R}$ ($\lpl\le \lmin$), which correspond
to a positive and a negative eigenfunction, respectively, and showed that
the positivity of these numbers is a sufficient condition for the solvability
of the related Dirichlet problem. Recently in \cite{QS} the second and the
third author  extended these results to arbitrary operators and studied the
properties of the eigenvalues and the eigenfunctions, in particular the
relation between the positivity of the eigenvalues and the validity of the
comparison principle and the Alexandrov-Bakelman-Pucci estimate, thus
obtaining extensions to nonlinear operators of the results of
Berestycki-Nirenberg-Varadhan in \cite{BNV}. In what follows we always assume that $F$ is indeed nonlinear in the sense that $\lpl< \lmin$
 --- note the results in \cite{QS} easily imply that $\lpl= \lmin$ can occur only  if all  linear operators
which appear in \re{formf} have the same principal eigenvalues {\it and}
eigenfunctions.

In the subsequent works \cite{S1}, \cite{FQS} we considered the Dirichlet problem
 \re{first} with $f$ independent of $u$, and we obtained a
number of results on the structure of its solution set, depending on the
position
of the parameter $\lambda$ with respect to the eigenvalues $\lpl$ and $\lmin$.
In particular,  we proved  that for each $\lambda$ in the closed interval $[\lpl,\lmin]$ and each $h\in L^p$, $p>N$, which is not a multiple of the first eigenfunction $\fipl$, there exists a critical number $t^*_{\lambda,F}(h)$ such that  the equation
\beq\label{ett}
F[u] + \lambda u = t\fipl+h\quad\mbox{in }\Omega\qquad u=0 \quad\mbox{on }\partial\Omega,
\eeq
has solutions for $t>t^*_{\lambda,F}(h)$ and has no solutions for $t<t^*_{\lambda,F}(h)$. We remark this is in sharp contrast with the case of linear $F$, say $F=\Delta$, when  \re{ett} has a solution if and only if $t=t^*_{\lambda,\Delta}(h)=-\int_\Omega (h\varphi_1)$ (we shall assume all eigenfunctions are normalized so that their $L^2$-norm is one).
Much more information on the solutions of \re{ett} can be found in \cite{S1} and \cite{FQS}.
The value of $t_{\lpl,F}^*(h)$ in terms
of $F$ and $h$ was computed by Armstrong
\cite{armstrong}, where he obtained an extension to HJB operators
of the Donsker-Varadhan minimax formula.

We now turn to the statements of our main results.
A standing assumption on the function $f$ will be the following
\begin{itemize}
\item[(F0)] $f:\bar\Omega\times\R\to\R$ is  continuous and  sub-linear in $u$ at infinity :
$$
\lim_{|s|\to\infty}\frac{f(x,s)}{s}=0 \quad\mbox{ uniformly in } \; x\in\bar\Omega.
$$
\end{itemize}

\bremark
For continuous  $f$ it is known (\cite{C}, \cite{swiech}, \cite{winter}) that all viscosity solutions of \re{first} are actually strong, that is, in $W^{2,p}(\Omega)$, for all $p<\infty$. Without serious additional complications we could assume that the dependence of $f$ in $x$ is only in $L^p$, for some $p>N$.
\eremark

\bremark
 Some of the statements below can be divided into subcases by supposing that $f$ is sub-linear in $u$ only as $u\to\infty$ or as $u\to-\infty$ (such results for the Laplacian can be found in \cite{CMN}, \cite{CdF}). We have chosen to keep our theorems as simple as possible.
\eremark

Now we introduce the hypotheses which
extend the Landesman-Lazer
assumptions  \re{ll} for the Laplacian to the case of general
HJB operators.
From now on we write
the critical $t$-values at resonance as
 $\tplus=t_+^*(h)=t_{\lpl,F}^*(h)$ and $\tminus=t_-^*(h)=t_{\lmin,F}^*(h)$, and $p>N$ is a fixed number. We assume there are
\begin{itemize}
\item[\hip] a function $c_+\in  L^p(\Omega)$, such that
$
c_+(x)\le f_+(x)$ in  $\Omega$  and $ \tplus(c_+)<0$.
\item[\him]
a function $c^-\in  L^p(\Omega)$, such that
 $c^-(x)\ge f^-(x)$   in  $ \Omega$  and
$\tminus(c^-)>0.$
\item[\hdp]
a function  $c^+\in  L^p(\Omega)$, such that
 $c^+(x)\ge f^+(x)$  in  $ \Omega$ and
$ \tplus(c^+)>0$.
\item[\hdm]
a function  $c_-\in  L^p(\Omega)$, such that
 $c_-(x)\le f_-(x)$ in  $ \Omega$ and $\tminus(c_-)<0$.
\end{itemize}
\bremark
Note that, decomposing  $h(x)= \left(\int_\Omega h\phiplus\right)\,\phiplus(x)+ h^\perp(x)$, where $\fipl$ is the eigenfunction associated to $\lambda_1^+$, we clearly have
 \begin{equation}\label{clear}t^*_\lambda(h) = t^*_\lambda(h^\perp) - \int_\Omega (h\phiplus)\quad\mbox{ for } \mbox{ each }\; \lambda\in [\lpl,\lmin].
  \end{equation}
  So when $F=\Delta$ hypotheses \hip-\him \ and \hdp-\hdm \ reduce to the classical Landesman-Lazer conditions \re{ll}, since for the Laplacian the critical $t$-value of a function orthogonal to $\varphi_1$ is always zero, by the Fredholm alternative.

We further observe that whenever
one of the limits $f_\pm, f^\pm$ is infinite, a
function $c_\pm, c^\pm$ with the required in \hip-\him, \hdp-\hdm \ property
always exists, while if any of $f_\pm, f^\pm$ is in $L^p(\Omega)$, we take the corresponding $c$ to be equal to this limit. Note also that the strict inequalities in  \hip-\him \ and \hdp-\hdm \ are important, see Section \ref{examples}.
\eremark

Throughout the paper we denote by  ${\cal S}$ the set of
all pairs
$(u,\lambda)\in C(\bar\Omega)\times \R$
which satisfy equation
\equ{first}. For any fixed $\lambda$ we set ${\cal S}(\lambda)= \{u\:|\: (u,\lambda)\in {\cal S}\}$ and if ${\cal C}\subset{\cal S}$ we denote $\calc(\lambda)= \calc \cap \cals(\lambda)$.

Our first result gives a statement of existence of solutions for $\lambda$ around $\lpl$ and $\lmin$, under the above Landesman-Lazer type hypotheses.
We recall that for some constant $\delta_0>0$
(all constants in the paper will be allowed to depend on
$N,\lambda, \Lambda,\gamma$, diam$(\Omega)$),  $\lpl,\lmin$ are
the only eigenvalues of $F$ in the interval $(-\infty,\lmin+\delta_0)$ --
see Theorem 1.3 in \cite{QS}.

\begin{teo}\label{teo1}
Assume (F0) and  \hip. Then  there exists  $\delta>0$ and
two  disjoint closed connected sets of solutions of \equ{first},
${\cal C}_1,\calc_2\subset {\cal S} $  such that
\begin{enumerate}
\item  ${\cal C}_1(\lambda)\not =\emptyset$ for all
$\lambda\in (-\infty, \lambdaplus]$,
\item
  ${\cal C}_1(\lambda)\not =\emptyset$  and
${\cal C}_2(\lambda)\not =\emptyset$ for all
$\lambda\in( \lambdaplus, \lambdaplus+\delta)$.
\end{enumerate}

The set $\calc_2$ is a branch of solutions "bifurcating from plus infinity to
the right of $\lambda_1^+$",
that is, ${\cal C}_2 \subset C(\overline \Omega)\times(\lambda_1^+,\infty)$ and
 there is a sequence  $\{(u_n,\lambda_n)\}\in \calc_2$
such that $\lambda_n\to \lpl$ and  $\|u_n\|_\infty\to \infty$. Moreover, for every
sequence $\{(u_n,\lambda_n)\}\in \calc_2$ such that $\lambda_n\to \lpl$ and  $\|u_n\|_\infty\to \infty$, $u_n$ is positive in $\Omega$, for $n$ large enough.

If we assume \him\
then
 there is a branch of solutions of \equ{first} "bifurcating from
minus infinity to the right of $\lmin$", that is, a connected set
${\cal C}_3\subset {\cal S} $ such that
 ${\cal C}_3\subset C(\overline \Omega)\times(\lambda_1^-,\infty)$
for which there is a sequence  $\{(u_n,\lambda_n)\}\in \calc_3$
such that $\lambda_n\to \lmin$ and  $\|u_n\|_\infty\to \infty$. Moreover, for every
sequence $\{(u_n,\lambda_n)\}\in \calc_3$ such that $\lambda_n\to \lmin$
and  $\|u_n\|_\infty\to \infty$, $u_n$ is negative in $\Omega$ for $n$ large.
%
\end{teo}

Under the sole hypothesis \hip\,
it cannot be guaranteed that the
sets of solutions
${\cal C}_1$ and ${\cal C}_2$ extend much beyond $\lpl$. This important fact
will be proved in Section
\ref{examples}, where we find $\delta_0>0$ such that for each $\delta\in(0,\delta_0)$ we can construct a nonlinearity
$f(x,u)$ which satisfies (F0), \hip\, and \him,
but for which   $\cals(\lambdaplus+\delta)$ is empty.

It is clearly important to give hypotheses on $f$ under which we can get a
global result, that is, existence of continua of solutions which extend
over the gap between the two principal eigenvalues
(this gap accounts for the nonlinear nature of the HJB operator !).
The next theorems deal with that question, and use the following additional
assumptions.
\begin{itemize}
\item[(F1)] $f(x,0)\ge 0$ and $f(x,0)\not\equiv0$ in $\Omega$;
\item[(F2)] $f(x,\cdot)$ is locally Lipschitz, that is, for each $R\in \mathbb{R}$ there is $C_R$ such that
$|f(x,s_1)-f(x,s_2)|\le C_R|s_1-s_2|$ for all $s_1,s_2\in  (-R,R)$ and $x\in \overline{\Omega}$.
\end{itemize}
A discussion on these hypotheses, together with examples and counterexamples,
will be given in Section \ref{examples}.
\begin{teo}\label{teo2}
Assume (F0), (F1), (F2),  \hip \, and \him \, hold. Then  there exist a
constant $\delta>0$ and
three  disjoint closed connected sets of solutions
 ${\cal C}_1, \calc_2$,  $\calc_3 \subset {\cal S}$,
 such that
\begin{enumerate}
\item[1.] $ {\cal C}_1(\lambda)\not =\emptyset$ for all $\lambda\in
(-\infty,  \lambdaplus]$,
\item[2.]
$ {\cal C}_i(\lambda)\not =\emptyset$ , $i=1,2$, for all $\lambda\in
( \lambdaplus, \lambdaminus]$,
\item[3.]
 $ {\cal C}_i(\lambda)\not =\emptyset$ , $i=1,2,3$,  for all $\lambda\in
( \lambdaminus, \lambdaminus+\delta)$.
\end{enumerate} The  sets $\calc_2$ end $\calc_3$ have the same "bifurcation from infinity"  properties as in the previous theorem.
\end{teo}
While Theorem \ref{teo2} deals with bifurcation branches going to the right
of the corresponding eigenvalues, the next theorem takes care of the
case where the branches go to the left of the eigenvalues.

\begin{teo}\label{teo3}
Assume (F0), (F1), (F2), \hdp\,  and \hdm \, hold. Then
there exist $\delta>0$ and disjoint closed connected sets of solutions
  ${\cal C}_1, \calc_2\subset {\cal S} $
 such that
\begin{enumerate}
\item[1.] ${\cal C}_1(\lambda)\not =\emptyset$ for all
$\lambda\in (-\infty,  \lambdaplus-\delta]$,
\item[2.]
 ${\cal C}_1(\lambda)\not =\emptyset$, $\calc_2(\lambda)$ contains at least
two elements
 for all $\lambda\in ( \lambdaplus-\delta, \lambdaplus)$,
and  $\calc_2$ is a branch "bifurcating from plus infinity to the left of
$\lambda_1^+$".
\item[3.]
${\cal C}_1(\lambda)\not =\emptyset$ and ${\cal C}_2(\lambda)\not =\emptyset$
 for all $\lambda\in [ \lambdaplus, \lambdaminus)$, and either:

(i)  ${\cal C}_1$ is the branch "bifurcating from minus infinity to the left of
$\lambda_1^-$"

(ii) There is a closed connected  set of solutions ${\cal C}_3\subset {\cal S}$,
 disjoint of
${\cal C}_1$  and ${\cal C}_2$,
"bifurcating from minus infinity to the left of $\lambda_1^-$"
  such that
 ${\cal C}_3(\lambda)$
has at least two elements for all $\lambda\in
(\lambda_1^--\delta,\lambda_1^-).$
\item[4.] ${\cal C}_2(\lambda)\not=\emptyset$  for all
$\lambda\in [\lambdaminus, \lambdaminus+\delta]$.
In case ii) in 3., ${\cal C}_2(\lambda)\not=\emptyset$ and
${\cal C}_3(\lambda)\not=\emptyset$ for all
$\lambda\in [\lambdaminus, \lambdaminus+\delta]$.
\end{enumerate}\end{teo}

Note alternative 3. (ii) in this theorem is somewhat anomalous. While we are able to exclude it in a number of particular cases (in particular for  the model nonlinearities which satisfy the hypotheses of the theorem), we do not believe it can be ruled out in general. See Proposition \ref{propsuppl} in Section \ref{sec_proofmain}.

Going back to the case when $F$ is linear, a well-known "rule of thumb"
states that the number of expected solutions of \re{first} changes by two
when the parameter $\lambda$ crosses the first eigenvalue of $F$.
An heuristic way of interpreting our theorems is that, when $F$ is a supremum of linear operators, crossing a "half"-eigenvalue leads to a change of the number of solutions by one.

The following graphs illustrate our theorems.

\begin{figure}[h]
\includegraphics[width=14cm, height=3cm]{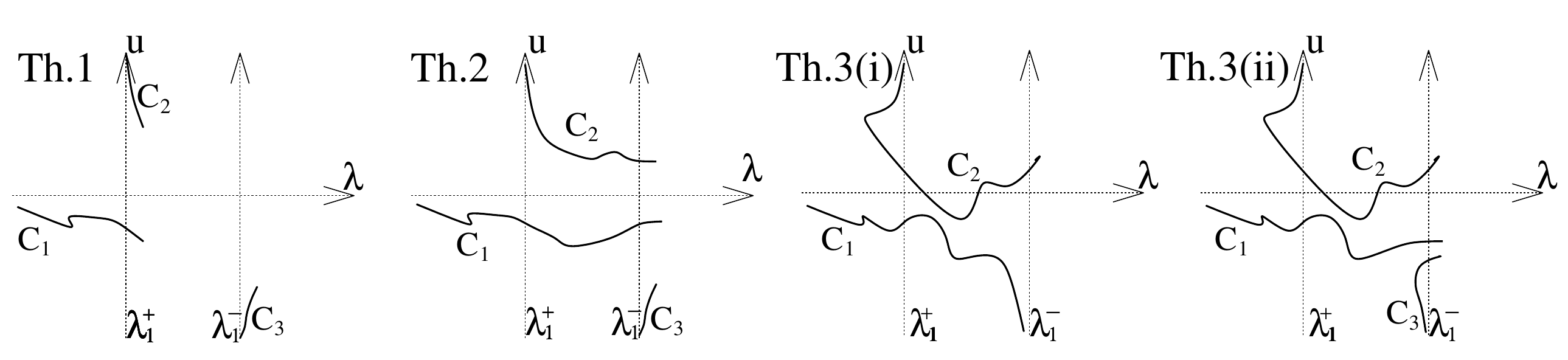}\label{fig1}
\end{figure}

The paper is organized as follows. The next section contains some definitions,
known results, and continuity properties of the critical values
 $t^*$. In Section \ref{sec_apriori} we obtain a priori bounds for
the solutions of \equ{first}, and construct  super-solutions or sub-solutions in  the different cases. In Section \ref{sec_bifurc} bifurcation from infinity for HJB operators is established through the classical method of Rabinowitz, while in Section \ref{sec_branch} we construct and study a bounded branch of solutions of \equ{first}. These results are put together in  Section \ref{sec_proofmain}, where we prove  our main theorems. Finally, a discussion on our hypotheses and some examples which highlight their role are given in Section \ref{examples}.

\setcounter{equation}{0}
\section{Preliminaries and continuity of $t^*$}\label{sect3}

First, we list the properties shared by HJB operators of our type. The function
$F:S_N\times \R^N\times\R\times\Omega\to\R$ satisfies (with $S,T\in S_N\times \R^N\times\R$)
\begin{itemize}
\item[(H0)]  $F$ is
positively homogeneous of order 1 : $
F(tS,x)=tF(S,x)
$ for  $t\ge 0$.

\item[(H1)] There exist $\lambda,\Lambda,\gamma>0$ such that for $S=(M,p,u), T= (N,q,v)$
\begin{eqnarray}
\mm(M-N)-\gamma(|p-q|+|u-v|)
\le F(S,x)-F(T,x)\nonumber\\
\le
\mp(M-N)+\gamma(|p-q|+|u-v|).\nonumber
\end{eqnarray}
\item[(H2)] The function
$F(M,0,0,x)$ is continuous in $S_N\times\overline{\Omega}$.
\item[(DF)] We have
$-F(T-S,x)\le F(S,x)-F(T,x)
\le F(S-T,x)$ for all $S,T$.
\end{itemize}
In (H1)  $\mm$ and $\mp$ denote the Pucci extremal operators, defined as $ \mp(M)=
 \mathop{\sup}_{A\in {\cal A}}\mbox{tr}(AM)$,
$\mm(M)=\mathop{\inf}_{A\in {\cal A}}\mbox{tr}(AM), $ where $
{\cal A} \subset {\cal S}_N$ denotes the set of matrices whose
eigenvalues lie in the interval $[ \lambda, \Lambda ]$, see for instance \cite{CC}. Note under
(H0)  assumption (DF) is equivalent to the convexity of $F$  in $S$
 -- see  Lemma 1.1 in \cite{QS}.
Hence for each $\phi,\psi\in W^{2,p}(\Omega)$ we have the inequalities $F[\phi+\psi]\le F[\phi] + F[\psi]$ and $F[\phi-\psi]\ge F[\phi] - F[\psi]$.

We recall the definition of the principal eigenvalues of $F$ from \cite{QS}
$$
\lambda_1^+(F,\Omega)=\sup\{\lambda\,|\, \Psi^+(F,\Omega,\lambda)\not =\emptyset\},\quad
\lambda_1^-(F,\Omega)=\sup\{\lambda\,|\, \Psi^-(F,\Omega,\lambda)\not =\emptyset\},$$
where
$
\Psi^\pm(F,\Omega,\lambda)= \{\psi\in C(\overline{\Omega})\;|\;
\pm(\F{\psi}+\lambda\psi)\le 0,\quad \pm\psi>0 \mbox{ in } \Omega \}.$
Many properties of the eigenvalues (simplicity, isolation, monotonicity and continuity with respect to the domain, relation with the maximum principle) are established in Theorems 1.1 -- 1.9 of  \cite{QS}. We shall repeatedly use these results. We shall also often refer to the statements on the solvability of the Dirichlet problem, given in \cite{QS} and \cite{FQS}.


We  recall the following Alexandrov-Bakelman-Pucci (ABP) and $C^{1,\alpha}$ estimates, see \cite{GT}, \cite{CCKS},  \cite{winter}.

\begin{teo} \label{estimate}
Suppose $F$ satisfies (H0), (H1), (H2),
 and $u$ is a solution of
$\F{u}+c u=f(x)$ in $\Omega$, with
$u = 0$ {on } $\partial\Omega$.
Then there
exists $\alpha \in (0,1)$ and  $C_0 > 0$ depending on $N, \lambda, \Lambda,
\gamma, c$ and $\Omega$ such that $u\in C^{1,\alpha}(\bar\Omega)$, and
$$
  \|u\|_{C^{1,\alpha}(\bar\Omega)}\le C_0(\|u\|_{L^{\infty}(\Omega)} +
\|f\|_{L^{p}(\Omega)}).
$$
Moreover, if one  chooses $c=-\gamma$ (so that by (H1) $F-\gamma$ is proper) then this equation has a
unique solution which satisfies
$
  \|u\|_{C^{1,\alpha}(\bar\Omega)}\le C_0\|f\|_{L^{\infty}(\Omega)}.
$ More precisely, any solution of $F[u]-\gamma u \ge f(x)$ satisfies
$\displaystyle\sup_\Omega u \le \sup_{\partial\Omega}u + C \|f\|_{L^N}$.
\end{teo}

For readers' convenience we state a version of  Hopf's Lemma (for viscosity solutions it was proved in \cite{bardidalio}).
\begin{teo}\label{bdl}
Let $\Omega\subset\R^N$ be a regular domain and let $\gamma>0$,
$\delta>0$. Assume $w\in C(\overline\Omega)$ is a viscosity solution of
$
{\cal M}^-_{\lambda,\Lambda}(D^2w)-\gamma|Dw|-\delta w\le 0$ in $\Omega$,
and $w\ge 0$ in $\Omega$. Then either $w\equiv 0$ in $\Omega$ or $w>0$ in $\Omega$
and at any point $x_0\in \partial \Omega$ at which  $w(x_0)=0$ we have
$
\limsup_{t\searrow 0}\frac{w(x_0+t\nu)-w(x_0)}{t}<0,
$
where $\nu$ is the interior normal to $\partial \Omega$ at $x_0$.
\end{teo}
The next theorem is a consequence of the compact embedding $C^{1,\alpha}(\Omega)\hookrightarrow C^{1}(\Omega)$,  Theorem \ref{estimate}, and the convergence properties of viscosity solutions (see Theorem 3.8 in \cite{CCKS}).
\begin{teo}\label{conv}  Let $\lambda_n\to \lambda$ in $\mathbb{R}$ and $f_n\to f$ in $L^p(\Omega)$. Suppose  $F$ satisfies $(H1)$ and $u_n$ is a viscosity solution of $F[u_n] +\lambda_n u_n= f_n$ in $\Omega$, $u_n=0$ on $\partial \Omega$. If $\{u_n\}$ is bounded in $L^\infty(\Omega)$ then a subsequence of $\{u_n\}$ converges in  $C^{1}(\overline{\Omega})$ to a function $u$, which solves $F[u] +\lambda u = f$ in $\Omega$,
$u=0$ on $\partial \Omega$. \end{teo}

As a simple consequence of this theorem, the homogeneity of $F$ and the simplicity of the eigenvalues we obtain the following proposition.
\begin{prop}\label{propconv}  Let $\lambda_n\to \lambda_1^\pm$ in $\mathbb{R}$ and $f_n$ be bounded in $L^p(\Omega)$. Suppose  $F$ satisfies $(H1)$ and $u_n$ is a viscosity solution of $F[u_n] +\lambda_n u_n= f_n$ in $\Omega$, $u_n=0$ on $\partial \Omega$. If $\{u_n\}$ is unbounded in $L^\infty(\Omega)$ then a subsequence of $\frac{u_n}{\norm{u_n}}$ converges in  $C^{1}(\overline{\Omega})$ to  $\varphi_1^\pm$. In particular, $u_n$ is positive (negative) for large $n$, and for each $K>0$ there is $N$ such that $|u_n|\ge K\fipl$ for $n\ge N$. \end{prop}

For shortness, from now on the zero boundary condition on
$\partial \Omega$ will be understood in all differential
(in)equalities we write, and $\| \cdot \|$ will refer to the $\linf$-norm.

We devote the remainder of this section to the definition and
some basic continuity properties of the critical $t$-values for \equ{ett}.
These numbers are crucial in the study of existence of solutions at resonance and
in the gap between the eigenvalues.
 For  each $\lambda\in[\lpl,\lmin]$
and each $d\in L^p$, which is not a multiple of the first
eigenfunction $\fipl$, the  number
$$t^*_\lambda(d)=\inf\{t\in \mathbb{R}\, |\,  F[u] + \lambda u = s\fipl+d\mbox{ has solutions for } s\ge t\}$$
is well-defined and finite. The non-resonant case $\lambda\in(\lpl,\lmin)$
was considered in  \cite{S1}, while
the resonant case $\lambda=\lpl$ and $\lambda=\lmin$  was studied in \cite{FQS}.

In what follows we prove the  continuity of $t^*_\lambda: L^p(\Omega)\to \R$
for any fixed $\lambda\in[\lambdaplus,\lambdaminus]$. Actually,
  in \cite{S1}  the continuity of this function is proved
 for all
$\lambda\in(\lambdaplus,\lambdaminus)$, so we only need to take care of the
resonant cases $\lambda=\lpl$ and $\lambda=\lmin$, that is, to study $\tplus$ and $\tminus$.
In doing so, it is convenient to use the following equivalent definitions of
 $\tplus$ and $\tminus$ (see \cite{FQS})
 \begin{equation}\label{equidef}\begin{array}{c}t^*_+(d)=\inf\{t\in
\mathbb{R}\, |\,\mbox{ for each }  s> t \mbox{ and  } \lambda_n\nearrow\lpl
\mbox{ there exists } u_n \\ \mbox{such that } F[u_n] + \lambda_n u_n =
s\fipl+d\mbox{ and } \|u_n\| \mbox{ is bounded }\}\end{array}
\end{equation}
and
\begin{equation}\label{equidef1}
 \begin{array}{c}t^*_-(d)=\inf\{t\in \mathbb{R}\, |\,\mbox{ for each }  s> t \mbox{ and  } \lambda_n\searrow\lmin \mbox{ there exists } u_n \\ \mbox{such that } F[u_n] + \lambda_n u_n = s\fipl+d\mbox{ and } \|u_n\| \mbox{ is bounded }\}.\end{array}
 \end{equation}
\begin{proposition}\label{continuoust}
The functions $\tplus,\tminus: L^p(\Omega)\to \R$ are continuous.
\end{proposition}
\noindent
{\bf Proof.}
 If we assume $\tplus$ is not continuous,
then there is  $d \in L^p(\Omega)$, $\varepsilon>0$ and a
sequence $d_n\to d$ in $L^p(\Omega)$
such that either
 $\tplus(d_n)\ge  \tplus(d)+3\varepsilon$ for all $n\in\N$ or
  $\tplus(d_n)\le  \tplus(d)-3\varepsilon$  for all $n\in\N$.

First we suppose that $\tplus(d_n)\ge  \tplus(d)+3\varepsilon$ for all $n\in\N$.
Then, for any sequence $\lambda_m\nearrow\lambdaplus$ we find solutions $u_n^m$
of the equation
$$
\F{u_n^m}+\lambda_m u_n^m = (\tplus(d_n)-2\varepsilon)\phiplus+d
\,\, \mbox{ in } \Omega,
$$
and the sequence $\{ \|u_n^m\|\}$ is bounded as $m\to\infty$, for each fixed $n$ -- see \re{equidef}. We also consider the
solutions $w_n$ of $
F[w_n]-\gamma w_n=d_n-d$.
By Theorem \ref{estimate} we know that  $w_n\to 0$ in
$C^1(\bar\Omega)$. Then by the structural hypotheses on $F$ (recall $\F{u+v}\le \F{u}+\F{v}$) we have
\begin{eqnarray*}
F[u_n^m+w_n]+\lambda_m (u_n^m+w_n)&\le& (\tplus(d_n)-2\varepsilon)\phiplus+d_n+
(\gamma+\lambda_m)w_n\\ &\le& (\tplus(d_n)-\varepsilon)\phiplus+d_n,
\end{eqnarray*}
where the last inequality holds if $n$ is large, independently of $m$. Fix one such $n$.
On the other hand we can take
solutions $z_n^m$ of
$$
\F{z_n^m} +\lambda_mz_n^m=(\tplus(d_n)-\varepsilon)\phiplus+d_n \quad(\ge F[u_n^m+w_n]+\lambda_m (u_n^m+w_n)).
$$
By \re{equidef} for any $n$ we have
$\| z_n^{m}\|\to \infty$ as $m\to\infty$.
By the comparison principle (valid by $\lambda_m<\lpl$ and Theorem 1.5 in \cite{QS})
 we obtain
$
z_n^{m}\le u_n^{m}+w_n$ in $\Omega$, hence $
z_n^{m}$ is bounded above as $m\to\infty$. Since $
z_n^{m}$ is bounded below, by Theorem  1.7 in \cite{QS} and $\lambda_m\le \lpl<\lmin$, we obtain a contradiction.

Assume now that $\tplus(d_n)\le  \tplus(d)-3\varepsilon$.
Let $u_n$ be a solution
of
$$
\F{u_n}+\lambdaplus u_n =(\tplus(d)-2\varepsilon)\phiplus+d_n
\quad \mbox{ in } \Omega,
$$
which exists since $\tplus(d_n)<\tplus(d)-2\varepsilon$ -- Theorem 1.2 in \cite{FQS}. Let $w_n$ be the solution of
$
\F{w_n}+cw_n = d-d_n$ in $\Omega,
$
with $w_n\to 0$ in  $C^1(\overline\Omega)$. Then there exists $n_0$ large enough so that $(\lambdaplus +\gamma)w_{n_0}<\varepsilon\varphi_1^+$, and consequently $u_{n_0} + w_{n_0}$ is a super-solution of
\begin{equation}\label{super}
\F{u} +\lambdaplus u=(\tplus(d)-\varepsilon)\phiplus+d
\quad \mbox{ in } \Omega.
\end{equation}
 Now, if
$w$ is the solution of
$
\F{w} -\gamma w = -d$
 in $\Omega,
$
by defining $v_k=k\phiminus-w$ we obtain
$$
\F{v_k} +\lambdaplus v_k\ge k(\lambdaplus-\lambdaminus)\phiminus
+d-(\lambdaplus+\gamma) w
>(\tplus(d)-{\varepsilon})\phiplus+d,
$$
for
$k$ large enough. By taking $k$ large we also have $v_k<u_{n_0} - w_{n_0}$ in $\Omega$, thus equation \equ{super} possesses ordered super- and sub-solutions. Consequently
it has a solution (by Perron's method -- see for instance Lemma 4.3 in \cite{QS}), a contradiction with the definition of $\tplus(d)$.
This completes the proof of the continuity of the function $\tplus$.

The rest of the proof is devoted to the analysis of continuity of
$\tminus$.
Assuming $\tminus$ is not continuous, there is $\varepsilon>0$ and a sequence $d_n\to d$ in $L^p(\Omega)$
such that either $\tminus(d_n)\ge  \tminus(d)+3\varepsilon$ or
$\tminus(d_n)\le  \tminus(d)-3\varepsilon$.
In the first case, let us consider a sequence $\lambda_m\searrow \lambdaminus,$
and a solution $v_m$ of the equation
$$
\F{v_m}+\lambda_mv_m=(\tminus(d)+\varepsilon)\phiplus +d
 \quad \mbox{ in } \Omega,
$$
We recall $v_m$ exists, by the results in \cite{armstrong} and \cite{FQS}.
We have shown in \cite{FQS} that $\tminus(d_n) \ge \tminus(d)+3\varepsilon>
\tminus(d)+\varepsilon$ implies that $v_m$ can be chosen to be  bounded as
$m\to \infty$ (see \re{equidef1}).  Let $w_n$ be the solution to
$
\F{w_n}-\gamma w_n=d_n-d$ in $\Omega,$
 as above. Then
$z_{n_0}^m=v_m+w_{n_0}$ satisfies for some large $n_0$
$$
\F{z_{n_0}^m}+ \lambda_m z_{n_0}^m\le(\tminus(d)+\varepsilon)\phiplus+(\lambda_m+\gamma) w_n+d_n
\le
(\tminus(d)+2\varepsilon)\phiplus+d_n,
$$
since again $w_n\to 0$ in $C^1(\overline\Omega)$.
On the other hand, we consider a solution of
$$
\F{u_n^m}+\lambda_mu_n^m=(\tminus(d)+2\varepsilon)\phiplus+d_n
\quad \mbox{ in } \Omega.
$$ As  $\tminus(d_n)\ge \tminus(d)+3\varepsilon>\tminus(d)+2\varepsilon$
for all  $n$, the sequence $u_n^m$ is not bounded (again by \re{equidef1}
and \cite{FQS}) and
$u_n^m/\|u_n^m\|_\infty\to\phiminus$ as $m\to\infty$, for each fixed $n$. Therefore for large $m$ the function
$\Psi=u_{n_0}^m-(v_m+w_{n_0})<0$ and
$
\F{\Psi}+\lambda_m\Psi\ge 0,
$
which is a contradiction with the definition of $\lambdaminus$, since $\lambda_m>\lambdaminus$.

Let us
assume now that for  $\varepsilon>0$ and the sequence $d_n\to d$ in
$L^p(\Omega)$
we have $t_-^*(d_n)\le t_-^*(d)-3\varepsilon$, for all $n$. Let
$\lambda_m\searrow \lambda_1^-$
and $v_m$ be a solution of the equation
$
\F{v_m}+\lambda_mv_m=(\tminus(d)-\varepsilon)\phiplus +d$
 in $\Omega$ (by \re{equidef1} $v_m$ is unbounded),
and let $w_n$ be the solution to
$
\F{w_n}-\gamma w_n=d-d_n$ in $\Omega.
$ Then,   $v_m/\|v_m\|_\infty\to \varphi_1^-$ and $w_n\to 0$ in
$C^1(\overline\Omega)$.
We take a solution $u_n^m$ to
$$
\F{u_n^m}+\lambda_m u_n^m =(\tminus(d)-2\varepsilon)\phiplus+d_n
\quad \mbox{ in } \Omega,
$$
and note that, since $\tminus(d_n)\le t^*_-(d)-3\varepsilon<t^*_-(d)-2\varepsilon,$ for any given $n$ there exists a constant
$c_n$ such that
$
\|u_n^m\|_\infty\le c_n,$ for all $m$.
Now, as above, we define $\Psi=v_m-(u_n^m+w_n)$, and  see that
\begin{eqnarray*}
\F{\Psi}+\lambda_m\Psi & \ge & (\tminus(d)-\varepsilon)\phiplus +d-
(\tminus(d)-2\varepsilon)\phiplus -d_n-
\F{w_n}-\lambda_mw_n\\ &  \ge &
\varepsilon\varphi_1^+-(\lambda_m+\gamma)w_n.
\end{eqnarray*}
We choose $n$ large enough to have $(\lambda_m+\gamma)w_n<
\varepsilon\varphi_1^+$ in $\Omega$. Then, keeping $n$ fixed, we can choose  $m$
large enough to have $\Psi<0$ in $\Omega$,  a contradiction with
the definition of $\lambda_1^-$, since $\lambda_m>\lambda_1^-$.\hfill
$\Box$

\medskip

 Finally we prove that  the function $t^*_\lambda(d)$
is also continuous in  $\lambda$ at the end points of the interval $[\lpl,\lmin]$, when $d$ is kept fixed. This fact will be needed in Section \ref{examples}.
\begin{proposition}\label{frehg} For every $d\in L^p(\Omega)$
$$
\lim_{\lambda\searrow\lambdaplus}t^*_\lambda(d)=\tplus(d)\quad\mbox{and}\quad
\lim_{\lambda\nearrow\lambdaminus}t^*_\lambda(d)=\tminus(d).
$$
\end{proposition}
\noindent
{\bf Proof.}
Let us assume that there is $\varepsilon>0$ and
 a sequence $\lambda_n\searrow\lambdaplus$ such that
$t^*_{\lambda_n}< \tplus-\varepsilon$ (since $d$ is fixed, we do not write it explicitly). Then  by the definition of  $t^*_{\lambda_n}$ there is a function
 $u_n$ satisfying
$$
\F{u_n}+\lambda_n u_n = (\tplus-\varepsilon)\phiplus +d\quad\mbox{ in } \Omega.
$$
Since $\lambda_n\searrow\lambdaplus$, $u_n$ cannot be bounded, for otherwise
we get a contradiction with the definition of $\tplus$ by finding a solution
with $\bar t<\tplus$ -- from Theorem~\ref{conv}. Then by Proposition \ref{propconv}
$u_n/\|u_n\|_\infty\to\phiplus$,  $u_n$
is positive for large $n$, and
$$\F{u_n}+\lambdaplus u_n=(\tplus-\varepsilon)\phiplus +d +
(\lambdaplus-\lambda_n)u_n
<(\tplus-\varepsilon)\phiplus+d,$$
that is $u_n$ is a super-solution.
On the other hand, for $t>\tplus$,
let $u$ be a solution of
$
\F{u}+\lambdaplus u = t\phiplus+d$,
 in  $\Omega,$
then $u$ is a sub-solution for this equation with $(\tplus-\varepsilon)
\phiplus+d$ as a right hand side. By taking $n$ large enough, we have
$u_n\ge u$, so that the equation
$$
\F{u}+\lambdaplus u=(\tplus-\varepsilon)\phiplus +d
\quad
 \mbox{ in } \Omega
$$
has a solution, a contradiction with the definition of $\tplus$.

Now we assume that there is $\varepsilon>0$ and
 a sequence $\lambda_n\searrow\lambdaplus$ such that
$t_n=t^*_{\lambda_n}> \tplus+2\varepsilon$.
Let $v$ be a solution
to
$$
\F{v}+\lambda^+_1 v=(\tplus+\varepsilon/2)\phiplus +d
\quad \mbox{ in } \Omega,
$$
then
$
\F{v}+\lambda_n v=(\tplus+\varepsilon)\phiplus +d -\varepsilon/2\phiplus+(\lambda_n-\lambdaplus)v.
$
Since $\tplus+\varepsilon<t_{\lambda_n}^*-\varepsilon/2$, by
choosing $n$ large we find
$
\F{v}+\lambda_n v<(t_{\lambda_n}^*-\varepsilon/2)\phiplus +d,
$
so that $v$ is a super-solution of
\begin{equation}\label{tyjg}
\F{u}+
\lambda_n u=(t_{\lambda_n}^*-\varepsilon/2)\phiplus +d.
\end{equation}
 Next we consider  a solution $u_K$ of $\F{u}+(\lpl +\nu)u = K$ (where we have set $\nu = (\lmin-\lpl)/2>0$), for each $K>0$. Such a solution exists by Theorem 1.9
in \cite{QS}, and it further satisfies  $u_K<0$ in $\Omega$ and
$\norm{u_K}_\infty\to \infty$ as $K\to \infty$, so $|u_K|\ge C(K)\fipl$, where $C(K)\to \infty$ as $K\to \infty$. Let $w$ be the (unique)
solution of $\F{w}-\gamma w= -d$ in $\Omega$. Since $\F{u_K-w}\ge \F{u_K} - \F{w}$,
we easily see that the function $u_K-w$ is a sub-solution of \re{tyjg}
and $u_K-w<v$, for large $K$. Then Perron's method leads again to
a contradiction with  the definition of $t_{\lambda_n}$.
This shows  $t_\lambda^*$ is right-continuous at
$\lambdaplus$.

Now we prove the second statement of Lemma \ref{frehg}. Assume
there is $\varepsilon>0$ and
a sequence $\lambda_n\nearrow\lambdaminus$ such that
$t^*_{\lambda_n}<\tminus -\varepsilon$.
Let $u_n$ be a solution to
$$
\F{u_n}+\lambda_n u_n = (\tminus-\varepsilon)\phiplus+d
 \quad \mbox{ in } \Omega,
$$
Since $\lambda_n\to\lambdaminus$, $u_n$ cannot be bounded (as before) and then
$u_n/\|u_n\|_\infty\to\phiminus$. Thus, for large $n$ we have $u_n<0$ and
$$
\F{u_n}+\lambdaminus u_n=(\tminus-\varepsilon)\phiplus+d+(
\lambdaminus-\lambda_n)u_n<(\tminus-\varepsilon)\phiplus+d,
$$
so that $u_k$ is a super-solution for some large (fixed) $k$. Consider now a sequence
$\widetilde\lambda_n\searrow\lambdaminus$ and let $v_n$ be the solution to
$$
\F{v_n}+\widetilde \lambda_n v_n = (\tminus-\varepsilon)\phiplus+d
\quad \mbox{ in } \Omega,
$$
whose existence was proved in \cite{armstrong} and \cite{FQS}. Then  $v_n$ cannot be
bounded, so
$v_n/\|v_n\|_\infty\to\phiminus$, and for large $n$ we have
$$
\F{v_n}+ \lambdaminus v_n=(\tminus-\varepsilon)\phiplus+d
 +
(\lambdaminus-\tilde \lambda_n)v_n>(\tminus-\varepsilon)\phiplus+d,
$$
that is $v_n$ is a sub-solution. For the already fixed $u_k$, we can find $n$ sufficiently large so that $u_k>v_n$, which implies that the equation
$$
\F{u}+\lambdaminus u=(\tminus-\varepsilon)\phiplus+d
 \quad \mbox{ in } \Omega,
$$
has a solution, a contradiction with the definition of $\tminus$.

Finally, assume that there is $\varepsilon>0$ and
a sequence $\lambda_n\nearrow\lambdaminus$ such that
$t^*_{\lambda_n}>t^*_{\lambda_n}-2\varepsilon>\tminus +\varepsilon$. By Theorem 1.4 in \cite{FQS} we can find a function $u$  which solves  the equation
$
\F{u}+\lambdaminus u=(\tminus+\varepsilon)\phiplus+d$ in $\Omega.
$
Then
$$
\F{u}+\lambda_n u<(t_{\lambda_n}^*-\varepsilon)\phiplus+d-\varepsilon\phiplus+
(\lambda_n-\lambdaminus)\phiplus
<(t_{\lambda_n}^*-\varepsilon)\phiplus+d,
$$
so that $u$ is a super-solution of $\F{u}+\lambda_n u=(t_{\lambda_n}^*-\varepsilon)\phiplus+d$, for some large fixed $n$. As we explained above, since $\lambda_n<\lmin$, by Theorem 1.9 in \cite{QS} we can construct an arbitrarily negative sub-solution of this problem, hence a solution as well,  contradicting the definition of $t_{\lambda_n}^*$.\hfill
$\Box$

\setcounter{equation}{0}
\section{Resonance and a priori bounds}\label{sec_apriori}

In this section we assume that the nonlinearity $f(x,s)$ satisfies the one-sided Landesman-Lazer
conditions at resonance, that is, one of  \hip, \him, \hdp \, and \hdm.
Under each of these conditions we analyze the existence of super-solutions,  sub-solutions
and a priori bounds when $\lambda$ is close to the eigenvalues $\lambda_1^+$ and $\lambda_1^-$.
This information will allow us  to obtain
 existence of solutions by using degree theory and bifurcation
arguments. In particular we will get branches bifurcating  from infinity which curve right or left depending on the a priori bounds obtained here.

We start with the existence of a super-solution  and a priori bounds at $\lambdaplus$, under hypothesis \hip.
\begin{proposition}\label{fi+}
Assume $f$ satisfies (F0) and \hip.
Then there exists a super-solution $z$ such that
$
\F{z}+\lambda z< f(x,z) $ in $
\Omega$, for all
$\lambda\in (-\infty,\lambdaplus]$.
Moreover, for each $\lambda_0<\lambdaplus$
there exist $R>0$  and a super-solution $z_0$ such that if $u$ is
a solution of \equ{first}  with
$\lambda\in[\lambda_0,\lambdaplus]$, then
$
\|u\|\le R$ and $u\le z_0$ in~$\Omega.
$
\end{proposition}
\noindent
{\bf Proof.}  We first replace
$c_+$ by a more appropriate function: we claim
that for each  $\varepsilon>0$ there exist
$R>0$ and a function $d\in L^p(\Omega)$ such that
$$
\|d-c_+\|_{L^p(\Omega)}\le \vae \qquad\mbox{and} \qquad
u\ge R\phiplus \;\mbox{ implies }\; f(x,u(x))\ge d(x)\; \mbox{ in }\Omega.
$$
In fact, setting $\sigma= \frac{\varepsilon}{2|\Omega|^{1/p}}$, we can find
$s_0$ such that $
f(x,s)\ge c_+(x)-\sigma$ in $\Omega$, for all $s\ge s_0$.
Let $\Omega^R = \{x\in \Omega\:|\: R\fipl(x) > s_0\}$ and  define the function $d_R$ as $d_R(x)=c_+(x)-\sigma$ if $x\in \Omega^R$, and
$d_R(x)=-M$ for $x\in \Omega\setminus \Omega^R$, where $M$ is such that
$
f(x,s)\ge -M,$ for all $s\in [0,s_0].$ It is then trivial to check that the claim holds for $d=d_R$, if $R$ is taken such that $|\Omega\setminus\Omega^R|<(\vae/2M)^p$.

Now,  by \hip\ and the continuity of $\tplus$ (Proposition
\ref{continuoust}) we can fix $\varepsilon$ so small that the function $d$ chosen above satisfies
\begin{equation}\label{di}
 \tplus(d)<0.
\end{equation}
 Let
$z_n$ be a solution to
$$
\F{z_n}+\lambdaplus z_n=t_n\phiplus+ d
 \quad \mbox{ in } \Omega,
$$
 where $t_n \to  \tplus(d)<0$, $t_n \ge  \tplus(d)$, is a sequence such that $z_n$ can be chosen to be unbounded --- such a choice of $t_n$ and $z_n$ is  possible thanks to  Theorem 1.2 in \cite{FQS}. Then
 $z_n/\|z_n\|\to \phiplus$, which implies that for large $n$
$$
\F{z_n}+\lambdaplus z_n< d\quad \mbox{ and }\quad z_n\ge R\phiplus,
$$
by \re{di}, where $R$ is as in the claim above. Thus $z_n$ is a strict
super-solution and, since  $z_n$ is
positive,
$
\F{z_n}+\lambda z_n< f(x,z_n),
$
for all $\lambda\in (-\infty, \lambdaplus]$. From now on we fix one such  $n_0$ and drop the
index, calling the super-solution $z$.

Suppose there exists  an unbounded sequence  $u_n$ of
solutions to
$$
\F{u_n}+\lambda_n u_n= f(x,u_n)
\quad \mbox{ in } \Omega,
$$
with $\lambda_n\in[\lambda_0,\lambdaplus]$ and
$\lambda_n\to\bar\lambda$. If $\bar \lambda<\lambdaplus$ then a contradiction follows
since $\lambdaplus$ is the first eigenvalue (divide the equation by $\|u_n\|$ and let $n\to \infty$). If $\bar \lambda=\lambdaplus$
then $u_n/\|u_n\|\to \phiplus$, so that for $n$ large we have  $u_n> z$ and
$u_n\ge R\phiplus$, consequently $f(x,u_n)\ge d(x)$ in $\Omega$. Thus, setting
$w=u_n-z$ we get, by $\lambda_n\le\lpl, w>0$,
$$
\F{w}+\lambdaplus w\ge \F{u_n} - \F{z} +\lpl (u_n-z)> f(x,u_n)-d\ge 0.
$$
Since $w>0$, Theorem 1.2 in \cite{QS} implies the existence of a constant
$k>0$ such that $w=k\phiplus$,  a contradiction with the last strict inequality.
Now that we have an a priori bound for the solutions,
we may choose an appropriate $n_0$ for the definition of $z_0=z_{n_0}$, which makes it larger than all solutions.\hfill
$\Box$

Next we state an analogous proposition on the existence of a sub-solution to our problem
at $\lambdaminus$ under hypothesis \him.
\begin{proposition}\label{fi-}
Assuming that $f$ satisfies (F0) and \him,
there
 exist
a strict sub-solution $z$ such that
$
\F{z}+\lambda z> f(x,z)$ in $
\Omega$
for all $\lambda\in
(-\infty,\lambdaminus]$.
Moreover, for each $\delta>0$ there
 exist $R>0$ and
a sub-solution $z$ such that if $u$ solves \equ{first} with
$\lambda\in   [\lambdaplus+\delta,\lambdaminus]$ then
$
\|u\|_\infty\le R$ and $u\ge z$ in $
\Omega.$
\end{proposition}
\noindent
{\bf Proof.} By using essentially the same proof as in Proposition \ref{fi+},
we  can find $R>0$ and a function
$d\in L^p(\Omega)$ such that
 $\tminus(d)>0$,
and $
u\le -R\phiplus$ implies $ f(x,u(x))\le d(x)$ in $\Omega$.
Consider a sequence $t_n \searrow  \tminus(d)$ and
 solutions $z_n$  to
\begin{equation}\label{desub0}
\F{z_n}+\lambdaminus z_n=t_n\phiplus+ d
 \quad \mbox{ in } \Omega,
\end{equation}
 chosen so that $z_n$ is unbounded and
$z_n/\|z_n\|_\infty\to \phiminus$ -- see  Theorem 1.4 in~\cite{FQS} . Hence for  $n$ large enough
\begin{equation}\label{desub}
\F{z_n}+\lambdaminus z_n> d\qquad\mbox{ and }\qquad z_n\le -R\phiplus.\end{equation}
Thus $z_n$ is a strict
sub-solution and, since $z_n$ is
negative for sufficiently large~$n$,
$
\F{z_n}+\lambda z_n> f(x,z_n),
$
for all $\lambda\in (-\infty, \lambdaminus]$. Fix one such $n_0$ and set $z=z_{n_0}$.

If $u_n$ is an unbounded sequence  of
solutions to
$
\F{u_n}+\lambda_n u_n= f(x,u_n),
$ in $\Omega,
$
with $\lambda_n\in[\lambdaplus+\delta,\lambdaminus]$ and
$\lambda_n\to\bar\lambda$ we obtain a contradiction like in the previous proposition. Namely, if $\bar \lambda\in
[\lambdaplus+\delta,\lambdaminus)$ then the conclusion follows
since there are no eigenvalues in this interval. If $\bar \lambda=\lambdaminus$
then $u_n/\|u_n\|\to \phiminus$, so that for $n$ large $u_n< z$ and
$u_n\le -R\phiplus$, hence $f(x,u_n)\le d(x)$, which leads to the contradiction $F[z-u_n] + \lmin(z-u_n)\ge 0$ and $z-u_n>0$. Then, given the a priori bound, we can choose $n_0$ such that $z_{n_0}$ is smaller than all solutions.  \hfill
$\Box$

\medskip


The next two propositions are devoted to proving a priori bounds under hypotheses \hdp \, and \hdm.
\begin{proposition}\label{fd+}
Under assumption (F0) and  \hdp\,
for each  $\delta>0$ the solutions to \equ{first}
with $\lambda\in [\lambdaplus, \lambdaminus-\delta]$
are a priori bounded.
\end{proposition}
\noindent
{\bf Proof.} As in the proof of Proposition \ref{fi+}, we may choose $R>0$
and a
function $d$ so that $\tplus(d)>0$, that is, $
\int_\Omega d\varphi_1^+< \tplus(d^\perp)$ (recall \re{clear}),
and whenever $u\ge R\phiplus$ then $f(x,u)\le d$.
Let $\tilde t$ be fixed such that
$
\int_\Omega d\varphi_1^+<\tilde t<\tplus(d^\perp).
$
If the proposition were not true, then
there would be  sequences $\lambda_n\searrow\lambdaplus$ and $u_n$ of solutions to
$
\F{u_n}+\lambda_n u_n= f(x,u_n)$,
such that $u_n$ is unbounded. Then  $u_n/\|u_n\|\to \phiplus$, in particular, $u_n$ is  positive for large $n$. Then
$$
\F{u_n}+\lambdaplus u_n\le f(x,u_n)\le d< \tilde t\phiplus+d^\perp.
$$
that is,  $u_n$ is a super-solution of $
\F{u_n}+\lambdaplus u_n= \tilde t \phiplus+d^\perp
$.
Next, take the solution $w$ of
$
\F{w}-\gamma w = -d^\perp$ in $\Omega,
$
where, as before, $\gamma$ is the constant from (H1), so that $F-\gamma$ is proper. For $\alpha>0$ we define $v=\alpha\phiminus-w$, then
$$
\F{v}+\lambdaplus v\ge \alpha(\lambdaplus-\lambdaminus)\phiminus-(\lpl+\gamma)w+d^\perp,
$$
exactly like in the proof of Proposition \ref{continuoust}. If we choose $\alpha$ large enough, we see that $v$ is a sub-solution for
$
\F{u_n}+\lambdaplus u_n=\tilde t\phiplus+d^\perp
$, and $v$ is smaller than the super-solution we constructed before. The existence  of a solution to this equation contradicts
the definition of $\tplus(d^\perp)$ and $\tilde t < \tplus(d^\perp)$.\hfill
$\Box$

\begin{proposition}\label{fd-}
Under assumption (F0) and \hdm\,
there exists $\delta>0$ such that the solutions to \equ{first}
with $\lambda\in [\lambdaminus, \lambdaminus+\delta]$
are a priori bounded.
\end{proposition}
\noindent
{\bf Proof.} We proceed like in the proof of the previous proposition. Now $
\int_\Omega d\varphi_1^+> \tminus(d^\perp)$,
and whenever $u\le -R\phiplus$ then $f(x,u)\ge d$.
If $\tilde t$ is such that
$
\int_\Omega d\varphi_1^+>\tilde t>\tplus(d^\perp),
$
and we assume there are sequences $\lambda_n\searrow\lambdaminus$ and $u_n$  of solutions to
$
\F{u_n}+\lambda_n u_n= f(x,u_n)$ in $\Omega$,
such that $u_n$ is unbounded, we get $u_n/\|u_n\|\to \phiminus$, consequently
$$
\F{u_n}+\lambdaminus u_n> \tilde t\phiplus +d^\perp.$$
On the other hand
if $z$ solves
$
\F{z}+\lambdaminus z= \tilde t\phiplus+d^\perp$ in $\Omega$ (such $z$ exists by Theorem 1.4 in \cite{FQS}),
then
$
\F{u_n-z}+\lambdaminus (u_n-z)>0,
$
and $u_n-z<0$ in $\Omega$, for large $n$. Thus, we may apply Theorem 1.4 in \cite{QS} to
obtain $k>0$ so that $u_n-z=k\phiminus$, a contradiction with the strict inequality. \hfill
$\Box$

\setcounter{equation}{0}
\section{Bifurcation from infinity at $\lambda_1^+$ and $\lambda_1^-$.}\label{sec_bifurc}

In this section we prove the existence of unbounded branches of solutions of
\equ{first}, bifurcating from
infinity at the eigenvalues $\lambdaplus$ and $\lambdaminus$. Then, thanks to
the a priori bounds obtained in
Section \S 4, for the two types of Landesman-Lazer conditions
(see Propositions \ref{fi+}-\ref{fd-}), we may determine to which side
of the eigenvalues these branches curve.

We recall that   $F(M,q,u,x)+cu$ is
decreasing in $u$ for any $c\le -\gamma $, in other words, $F+c$ is a proper operator.
Given $v\in C^1(\bar\Omega)$ we consider the problem
\begin{equation}\label{equ01v}
 \F{u}+cu =(c-\lambda)v+f(x,v)
 \quad\mbox{in }\; \Omega,\qquad
u = 0 \quad\mbox{on }\; \partial\Omega,
\end{equation}
see Theorem \ref{estimate}. We define the operator $K: \R\times C^1(\bar\Omega)\to C^1(\bar\Omega)$ as follows: $K(\lambda, v)$ is the
 unique solution $u\in C^{1,\alpha}(\bar\Omega)$ of  \equ{equ01v}.
The operator $K$ is compact in view of Theorem \ref{estimate} and the compact
embedding $C^{1,\alpha}(\bar\Omega)\to C^1(\bar\Omega)$.
With these definitions, our equation \equ{first} is transformed into the
fixed point problem
$
u=K(\lambda,u)$, $u\in C^1(\bar\Omega),
$
with $\lambda\in \R$ as a parameter.
We are going to show that the sub-linearity of the function $f(x,\cdot)$, given by assumption (F0), implies
bifurcation at infinity at the eigenvalues of $F$. The proof  follows
the standard
procedure for the linear case, see for example \cite{rabinowitz} or
\cite{reichel}, so we shall be sketchy, discussing only the main  differences. We define
$$
G(\lambda,v)=\|v\|^2_{C^1}K(\lambda, \frac{v}{\|v\|^2_{C^1}}),
$$
for $v\not =0$, and $G(\lambda,0)=0$. Finding  $u\not =0$ such that $
u=K(\lambda,u)$ is
equivalent to solving
the fixed point problem
$v=G(\lambda,v)$, $v\in C^1(\bar\Omega)$,
for $v=u/\|u\|^2_{C^1}$. The important observation is that bifurcation from zero in $v$ is equivalent to bifurcation from infinity for $u$.

Let $u=G_0(\lambda,v)$
be the  solution of
the problem
\begin{equation}\label{equ01ve}
 \F{u}+cu =(c-\lambda) v
 \quad\mbox{in }\; \Omega,\qquad
u = 0 \quad\mbox{on } \partial\Omega,
\end{equation}
and set $G_1=G-G_0$, so that
$
G(\lambda,v)=G_0(\lambda,v)+G_1(\lambda,v).
$
\begin{lema}\label{supbif}
Under hypothesis (F0) we have
$
\displaystyle\lim_{\|v\|_{C^1}\to 0}\frac{G_1(\lambda,v)}{\|v\|_{C^1}}=0.
$
\end{lema}
\noindent
{\bf Proof.} Let $g=G(\lambda,v)$ and $g_0=G_0(\lambda,v)$. Then we have
$$
\frac{1}{\|v\|_{C^1}}\left(\F{g}-\F{g_0}+c(g-g_0)\right)=
\|v\|_{C^1}f(x,\frac{v}{\|v\|_{C^1}^2}).
$$
The right hand side here goes to zero as
$\|v\|_{C^1}\to 0$, by (F0).
Then by (DF)
$$
 \frac{1}{\|v\|_{C^1}} \left(\F{|g-g_0|+c|g-g_0|}\right)\ge -  \|v\|_{C^1}\left|f(x,\frac{v}{\|v\|_{C^1}^2})\right|,
$$so the ABP inequality (Theorem \ref{estimate}) implies
$$
\sup_\Omega \{\frac{1}{\|v\|_{C^1}}|g-g_0|\}\le C \|v\|_{C^1}  \|f(x,\frac{v}{\|v\|^2_{C^1}}) \|_{L^p},
$$
and  the result follows.\hfill
$\Box$

\medskip

The next proposition deals with  the equation $
v=G_0(\lambda,v)$, $v\in C^1(\bar\Omega)
$ (recall we want to solve $v=G_0(\lambda,v)+G_1(\lambda,v)$),
which is equivalent to
\begin{equation} F(D^2v,Dv,v,x) =-\lambda v
\,\,\mbox{in}\,\, \Omega,\label{equ01vev}\quad
v=0\,\, \mbox{on}\,\,\partial\Omega.
\end{equation}
\begin{proposition}\label{deg} There exists $\delta>0$
 such that for all $r>0$ and all
$\lambda\in (-\infty,\lambda_1^-+\delta)\setminus\{\lambda_1^+,\lambda_1^-\}$, the Leray-Schauder degree
 ${\rm deg}(I-G_0(\lambda,\cdot),B_r,0)$ is well defined. Moreover
$$
{\rm deg}(I-G_0(\lambda,\cdot),B_r,0) = \left\{ \begin{array}{ccl}
                    1 & \mbox{ if }&  \lambda <\lambda_1^+  \\
                   0 &   \mbox{ if } &  \lambda_1^+<\lambda<\lambda_1^-\\
                 -1 &  \mbox{ if } & \lambda_1^-<\lambda<\lambda_1^-+\delta.
                 \end{array}
        \right.
$$
\end{proposition}
\noindent
{\bf Proof.}
We recall  it was proved in
\cite{QS} that problem \equ{equ01vev}
has only the zero solution in $(-\infty,\lambda_1^-+\delta)
\setminus\{\lpl,\lmin\}$, for  certain $\delta>0$. The compactness of $G_0$
follows from
Theorem \ref{estimate}, so the degree is well defined in the given ranges for $\lambda$.

Suppose $\lambda<\lambda_1^+$ and consider
the operator $I-tG_0(\lambda,\cdot)$ for $t\in [0,1]$. Since
 $t\lambda$ is not an eigenvalue of
\equ{equ01vev}, we have for  $t\in [0,1]$
$$
{\rm deg}(I-G_0(\lambda,\cdot),B_r,0)=
{\rm deg}(I-tG_0(\lambda,\cdot),B_r,0)=
{\rm deg}(I,B_r,0)=1.
$$

The case  $\lambda_1^+< \lambda<\lambda_1^-$ was studied in \cite{S1}. Consider the problem
\begin{equation}\label{equ01vevt}
 \F{u}+cu =(c-\lambda) v-t\varphi_1^+
\; \mbox{ in } \Omega,\qquad
u = 0 \;\mbox{ on } \partial\Omega
\end{equation}
for $t\in[0,\infty)$, whose unique solution is denoted by $\tilde{G}_0(\lambda,v,t)$.
It follows from the results  in \cite{QS}, \cite{S1} that  for $t>0$ the equation
\begin{equation}\label{equ01vevtu}
\F{u}+cu =(c-\lambda) u -t\varphi_1^+
\; \mbox{ in } \Omega,\qquad
u = 0 \;\mbox{ on } \partial\Omega
\end{equation}
does not have a solution. On the other hand, since $\lambda$ is not an
eigenvalue, there is $R>0$ such that the solutions of
\equ{equ01vevtu}, for $t\in [0,\bar t]$, are a priori bounded, consequently
\begin{eqnarray*}
{\rm deg}(I-G_0(\lambda,\cdot),B_r,0) &=&
{\rm deg}(I-\tilde{G}_0(\lambda,\cdot,0),B_R,0)\\
&=&
{\rm deg}(I-\tilde{G}_0(\lambda,\cdot,\bar t),B_R,0)=0.
\end{eqnarray*}

If $\lambda_1^-< \lambda<\lambda_1^-+\delta$ we proceed as in \cite{FQS},
 where the computation of the degree was done by making a homotopy with the Laplacian (see the proof of Lemma 4.2 in that paper).\hfill
$\Box$

Now we are in  position to apply the general theory of bifurcation to $v=G(\lambda,v)$, see for instance the surveys \cite{rabinowitz} and \cite{reichel},  and obtain
bifurcation branches  emanating from $(\lambda_1^+,0)$ and
$(\lambda_1^-,0)$, exactly like in \cite{BEQ}. In short, from $(\lambda_1^+,0)$  bifurcates a continuum of solutions of $v=G(\lambda,v)$, which is either unbounded in
$\lambda$, or unbounded in $u$, or connects to $(\overline \lambda, 0)$,
where $\overline \lambda\not =\lpl$ is an eigenvalue
(recall $\lpl$ and $\lmin$ are the only eigenvalues in
$(-\infty, \lmin + \delta)$, for some $\delta>0$). A similar situation occurs at
 $(\lmin,0)$.
Inverting the variables we obtain bifurcation at
infinity for our problem \equ{first}:
\begin{teo}\label{p1}
Under the hypotheses of Theorem \ref{teo1} there are two connected sets
${\cal C}_2$,  ${\cal C}_3\subset {\cal S}$ such that

1) There is a sequence $(\lambda_n, u_n)$ with $u_n\in
{\cal C}_2(\lambda_n)$ ($u_n\in {\cal C}_3(\lambda_n)$), and
 $\|u_n\|_\infty\to \infty$, $\lambda_n\to \lambda_1^+$ ($ \lambda_1^-$).

2) If
$(\lambda_n, u_n)$ is a sequence such that  $u_n\in {\cal C}_2(\lambda_n)$
(${\cal C}_3(\lambda_n)$),
 $\|u_n\|_\infty\to \infty$ and $\lambda_n\to \lambda_1^+$ ($ \lambda_1^-$), then $u_n$ is positive (negative) for large $n$.

3) The branch ${\cal C}_2$ satisfies one of the following alternatives, for some $\delta>0$:
(i) ${\cal C}_2(\lambda) \not = \emptyset$ for all  $\lambda\in (\lambda_1^+, \lambda_1^-+\delta)$ ;
(ii) There is $\lambda \in (-\infty, \lambda_1^-+\delta
]$ such that
$ 0\in {\cal C}_2(\lambda)$ ;
(iii) ${\cal C}_2(\lambda) \not = \emptyset$ for all
$\lambda\in (-\infty,\lambda_1^+)$ ;
(iv)
There is a sequence $(\lambda_n, u_n)$ such that $u_n\in {\cal C}_2(\lambda_n)$,
 $\|u_n\|_\infty\to \infty$, $\lambda_n\to \lambda_1^-$, and $\lambda_n\le \lambda_1^-$.

4)  The branch ${\cal C}_3$ satisfies one of the following
alternatives, for some
 $\delta>0$:
(i) ${\cal C}_3(\lambda) \not = \emptyset$ for all  $\lambda\in (\lambda_1^-,\lambda_1^-+\delta)$ ;
(ii) There is $\lambda \in (-\infty, \lambda_1^-+\delta
]$ such that
$ 0\in {\cal C}_3(\lambda)$ ;
(iii) ${\cal C}_3(\lambda) \not = \emptyset$ for all
$\lambda\in (-\infty,\lambda_1^-)$ ;
(iv)
There is a sequence $(\lambda_n, u_n)$ such that $u_n\in {\cal C}_3(\lambda_n)$,
 $\|u_n\|_\infty\to \infty$, and $\lambda_n\to \lambda_1^+$.
\end{teo}

We remark that (F1) excludes alternatives 3) (ii) and 4) (ii) in this theorem.

\setcounter{equation}{0}
\section{ A bounded branch of solutions }\label{sec_branch}

In this section we prepare for the proof of our main theorems by
establishing the existence of a continuum of solutions of \equ{first} which is not empty for all $\lambda\in(-\infty,\lpl+\delta)$, for some $\delta>0$.
Our first proposition concerns the behavior of solutions of
\equ{first} when $\lambda\to -\infty$.
\begin{proposition}\label{mininf}
1. Assume $f$ satisfies (F0).  Then there exists a constant $C_0>0$,
depending only on $F,f$, and $\Omega$, such that any solution of \equ{first}
satisfies $\|u\|_\infty\le C_0\lambda^{-1}$ as $\lambda\to -\infty$.

2. If in addition $f$ is Lipschitz at zero, that is, for some $\varepsilon>0$ and some $\overline{C}>0$ we have $|f(x,s_1)-f(x,s_2)|\le \overline{C}|s_1-s_2|$ for $s_1,s_2\in(-\varepsilon,\varepsilon)$, then \equ{first} has at most one solution when $\lambda$ is sufficiently large and negative. \end{proposition}

\noindent{\bf Proof.} 1. Let $u_\lambda$ be a sequence of solutions of \equ{first}, with $\lambda\to-\infty$. We first claim that $\|u_\lambda\|_\infty$ is bounded. Suppose this is not so, and say $\|u_\lambda^+\|_\infty\to\infty$ (with the usual notation for the positive part of $u$). Then, setting $v_\lambda= u_\lambda/\|u_\lambda^+\|_\infty$,  on the set $\Omega_\lambda^+=\{u_\lambda>0\}$ we have the inequality
$$
F[v_\lambda]-\gamma v_\lambda \ge \frac{f(x,u_\lambda)}{\|u_\lambda^+\|_\infty}\to 0,\; \mbox{ as } \lambda\to -\infty.
$$
The ABP estimate (see Theorem \ref{estimate}) then implies $\sup_{\Omega_\lambda^+} v_\lambda \to 0$, which is a contradiction with $\sup_{\Omega_\lambda^+} v_\lambda =1$. In an analogous way we conclude that $\|u_\lambda^-\|_\infty$ is bounded.

Hence there exists a constant $C$ such that $|f(x,u_\lambda(x))|\le C$ in $\overline{\Omega}$, so
$$
F[u_\lambda] - \gamma u_\lambda  \ge -(\lambda + \gamma)u_\lambda - C \ge 0\qquad \mbox{on the set }\; \tilde\Omega_\lambda,
$$
where $\tilde\Omega_\lambda = \{ u_\lambda>C/(|\lambda|+\gamma)\}$. Applying the maximum principle or the ABP inequality in this set implies it is empty, which means $u_\lambda\le C/(|\lambda|+\gamma)$ in $\Omega$. By the same argument we show $u_\lambda$ is bounded below, and  1. follows.

2. From statement 1. we conclude that for $\lambda$ small,
 all solutions of \equ{first} are in $(-\varepsilon,\varepsilon)$. If $u_1,u_2$ are two solutions of \equ{first} then for $|\lambda|>\gamma + \overline{C}$ we have $F[u_1-u_2] -\gamma (u_1-u_2)\ge 0$ on $\{u_1>u_2\}$ which means this set is empty. \hfill $\Box$
\medskip

The next result is stated in the framework of Theorem \ref{teo1} and  gives a bounded family of solutions $(u_\lambda,\lambda)$, for
$\lambda\in(-\infty,\lambda_1^++\delta)$. No assumption of Lipschitz continuity on $f$ is needed.
\begin{proposition}\label{degree}
Assume  $f$ satisfies
(F0) and ($F^\ell_+$). Then there is a connected subset ${\cal C}_1$ of ${\cal S}$ such that
${\cal C}_1(\lambda)\not =\emptyset$, for all
$\lambda\in (-\infty, \lambda_1^++\delta)$.
\end{proposition}
\noindent
{\bf Proof.}
According to Proposition \ref{fi+}, given $\lambda_0<\lambdaplus$, there is $R>0$ so that
all solutions of \equ{first} with $\lambda\in [\lambda_0, \lambda_1^+]$
belong to the ball $B_R$. In particular,
the equation does not have a solution  $(\lambda, u) \in [\lambda_0, \lambda_1^+]\times\partial B_R$. Moreover, there is
$\delta>0$ such that \equ{first} does not have a solution in
$ [\lambda_1^+, \lambda_1^++\delta]\times\partial B_R$ --- otherwise we obtain a contradiction by a simple passage to the limit.
Consequently the degree
${\rm deg}(I-K(\lambda,\cdot), B_R,0)$ is well defined
for all  $\lambda\in
[\lambda_0, \lambda_1^++\delta]$ ($K(\lambda,\cdot)$ is defined in the previous section). We claim that its value is $1$.

 To compute this degree, we fix
$\lambda< \lambda_1^+$ and analyze the equation
$$\F{u}+\lambda u=sf(x,u)\qquad \mbox{in }\; \Omega,$$
for $s\in [0,1]$.
Since $\lambda$ is not an eigenvalue of $F$ in $\Omega$,
 the solutions of this equation are a
priori bounded, uniformly in $s\in [0,1]$, that is, there is $R_1\ge R$, such that no solution of the equation exists outside of the open ball $B_{R_1}$.
Given $v\in C^1(\overline \Omega)$ we denote by $K_s(\lambda,v)$ the
unique solution of the equation
$\F{u}+\lambda u=sf(x,v) $ in $\Omega$.
Then we have
\begin{eqnarray*}
{\rm deg}(I-K(\lambda,\cdot), B_R,0) &=&
{\rm deg}(I-K_1(\lambda,\cdot), B_{R_1},0)\\
&=&{\rm deg}(I-K_0(\lambda,\cdot), B_{R_1},0)=1,
\end{eqnarray*}
where the last equality is given by Proposition \ref{deg}. Hence, again by the homotopy invariance of the degree, we have ${\rm deg}(I-K(\lambda,\cdot), B_R,0)=1$, for all $\lambda\in (\lambda_0,\lpl+\delta)$.

The last fact together with standard degree theory implies that for every $\lambda\in
[\lambda_0, \lambda_1^++\delta]$ there is at least one $(\lambda,u)$,
solution of \equ{first}, and,  moreover,
there is a connected subset ${\cal C}_1$ of ${\cal S}$ such that
${\cal C}_1(\lambda)\not =\emptyset$ for all $\lambda$
 in the interval
$[\lambda_0, \lambda_1^++\delta]$. Since $\lambda_0$ is arbitrary,
we can use the same argument for each element of a sequence  $\{\lambda_0^n\}$,
with $\lambda_0^n\to -\infty$. Then,
by a limit argument (like the one in the proof of Theorem 1.5.1 in \cite{FQS}), we find a connected set ${\cal C}_1$
with the desired properties.\hfill
$\Box$

\medskip

Next we study a branch of solutions driven by a family
of super- and sub-solutions, assuming that $f$ is locally Lipschitz continuous. In this case the statement of the previous proposition can be made more precise.
Specifically,  we assume that $f$ satisfies  (F2),  and there exist
 $\underline u, \overline u \in C^1(\overline \Omega)$,
such that $\overline u$ is a super-solution  and $\underline u$ is
 a sub-solution of  \equ{first}, for all $\lambda\le \bar\lambda$,
where $ \bar\lambda$ is fixed. We further
assume that $\underline u$ and  $\overline u$  are  not solutions of
\equ{first}, and
\begin{equation}\label{ss1}
\underline u< \overline u \mbox{ in }  \Omega,\qquad
\underline u= \overline u=0 \,\, \mbox{ and }  \,\,
\frac{\partial \underline u}{\partial \nu}<
\frac{\partial \overline u}{\partial \nu} \,\, \mbox{ on } \,\, \partial\Omega.
\end{equation}
We define the set
\begin{equation}\label{o}
{\cal O}=\{ v\in C^1(\overline \Omega)\;|\; \underline u< v<\overline u\,\,\mbox{ in }  \Omega\,\,
\mbox{ and } \,\, \frac{\partial \underline u}{\partial \nu}<
\frac{\partial v}{\partial \nu} <\frac{\partial \overline u}{\partial \nu} \,\,
\mbox{on} \,\, \partial \Omega\},
\end{equation}
which is open in $C^1(\overline \Omega)$.
Since $
{\cal O}$ is bounded in $C(\overline \Omega)$, we see that for every
$\lambda_0<\bar\lambda$ the set of solutions of \equ{first} in
$[\lambda_0, \bar\lambda]\times {\cal O}$ is bounded in
$C^1(\overline \Omega)$, that is,  all solutions
of  \equ{first} in
$[\lambda_0, \bar\lambda]\times {\cal O}$ are inside the ball
$B_R$, for some $R>0$.

\begin{lema}\label{supsubdeg}
With the definitions given above,
we have
$$
{\rm deg}(I-K(\lambda,\cdot), {\cal O}\cap B_R,0)=1, \qquad \mbox{ for all } \lambda\in
[\lambda_0, \bar\lambda].
$$
\end{lema}
\noindent
{\bf Proof.}
First we have to prove that the degree is well defined. We just need to show
that there are no fixed points of $K(\lambda,\cdot)$ on the boundary
of ${\cal O}\cap B_R$.
For this purpose it is enough to prove that, given $v\in C^1(\overline \Omega)$
such that $\underline u\le v\le \overline u$ in $\Omega$, we have
$\underline u < K(\lambda, v)< \overline u$ in $\Omega$. In what follows we write
$u= K(\lambda, v)$.

By (F2) we can assume that the negative number $c$, chosen
in Section \ref{sec_bifurc}, is such that the function $s\to f(x,s)+(c-\lambda)s$
is decreasing, for $s\in (-\tau,\tau)$, where $\tau=\max\{\|\underline{u}\|_\infty,\|\overline{u}\|_\infty\}$. Then
 \begin{eqnarray*}
\F{u}&=&\F{u}-f(x,\overline u)-
(c-\lambda) \overline u +f(x,\overline u)+(c-\lambda)\overline u\\
&\ge & \F{u}-f(x,v)-
(c-\lambda) v +f(x,\overline u)+(c-\lambda)\overline u\\
&=& -cu+f(x,\overline u)+(c-\lambda)\overline u\\
&=&c(\overline u-u)+f(x,\overline u)-\lambda\overline u\ge  \F{\overline u}+ c(\overline u- u).
\end{eqnarray*}
By (H1)  this implies
$
{\cal M}^+(D^2(u-\overline u))+\gamma |Du-D\overline u|+(\gamma-c)(u-\overline u)
>0
$
in $\Omega$. It follows from  Theorem \ref{bdl} that $u<\overline u$ in  $\Omega$ and $
\frac{\partial u}{\partial \nu} < \frac{\partial \overline u}{\partial \nu}$ on $\partial\Omega$.
The other inequality is obtained similarly.

By using its homotopy invariance, the degree we want to compute is equal to the degree at $\lambda_0$. But the latter was shown to be one in
the proof of Proposition \ref{degree}, which completes the proof of the lemma.\hfill
$\Box$

\medskip

Now we can state a proposition on the existence of a branch of solutions for
$\lambda \in (-\infty, \bar \lambda]$, whose proof is a direct consequence of
Lemma \ref{supsubdeg} and general degree arguments.
\begin{proposition}\label{boundedbranch}
Assume  $f$ satisfies
(F0) and (F2). Suppose there are functions
 $\underline u, \overline u \in C^1(\overline \Omega)$
such that $\overline u$ is a super-solution  and $\underline u$ is
 a sub-solution of  \equ{first} for all $\lambda\le \bar\lambda$, these functions are  not solutions of \equ{first} and
 satisfy \equ{ss1}.
Then there is a connected subset ${\cal C}_1$ of ${\cal S}$ such that
${\cal C}_1(\lambda)\not=\emptyset$ for all $\lambda\in (-\infty, \bar\lambda)$
and   each $u\in {\cal C}_1(\lambda)$ is such that $\underline u\le u \le \overline u$.
\end{proposition}
\bremark \label{rem1} In the next section we  use this proposition with  appropriately chosen  sub-solutions and super-solutions.
\eremark
\bremark \label{rem2}
The branch ${\cal C}_1$ is isolated of other branches of
solutions by the open set ${\cal O}$,
since we know there are no solutions on $\partial {\cal O}$.
\eremark

\setcounter{equation}{0}
\section{Proof of the main theorems}\label{sec_proofmain}

In this section  we put together  the bifurcation branches emanating from
infinity obtained in Theorem \ref{p1} with the bounded branches constructed in Section \ref{sec_branch}, and study their properties.
\medskip

\noindent
{\bf Proof of Theorem \ref{teo1}.}
This theorem is a  consequence of Proposition \ref{degree},
for the definition
of ${\cal C}_1$, and of Theorem \ref{p1} 1)-2), for the definition of ${\cal C}_2$ and ${\cal C}_3$.  Both ${\cal C}_2$ and ${\cal C}_3$
curve to the
right of $\lambda_1^+$ and  $\lambda_1^-$, respectively --  as a consequence of the a priori bounds
obtained in Proposition \ref{fi+} and \ref{fi-}.
\hfill $\Box$

\medskip

\noindent
{\bf Proof of Theorem \ref{teo2}.}
We first  construct  the branch ${\cal C}_1$, through
Proposition \ref{boundedbranch}.
In view of (F1)
we may take as super-solution   the function
$\overline u\equiv 0$.
In order to define the corresponding sub-solution we use Proposition \ref{fi-},
where a sub-solution is constructed for all $\lambda\in (-\infty,\lambda^-_1].$
We can rewrite  inequality \equ{desub0} in the following way
$$
\F{z_n}+(\lambdaminus +\delta)z_n= \tilde t_n\phiplus+ d,
\quad \mbox{ with }\quad
\tilde t_n(x)=\frac{\delta z_n(x)}{\phiplus(x)}+t_n.
$$
Since $z_n/\|z_n\|_\infty\to \phiminus<0$ in $C^1(\bar\Omega)$ we find that for some $c>0$
$$
\frac{|z_n(x)|}{\|z_n\|_\infty \phiplus(x)}\le c,\quad \forall x\in\Omega.
$$
Consequently, once $n$ is chosen so that  \equ{desub} holds, we can fix
$\delta>0$ such that
$\tilde t_n(x)\ge -\delta c + \tminus(d)>0$, for all  $x\in\Omega$,
which means that $z_n$ is a sub-solution also for $F[u] +(\lambda_1^-+\delta)u = f(x,u)$, as in the proof of Proposition \ref{fi-}.

Now we define  $\underline u=z_n$, chosen as above, and take $\bar \lambda=
\lambda_1^-+\delta$ in Proposition~\ref{boundedbranch}.
Clearly $\underline u$ and
$\overline u$ satisfy also \equ{ss1}, so the
 existence of the branch ${\cal C}_1$ (with the properties stated in
Theorem \ref{teo2}) follows from  Proposition
\ref{boundedbranch}.

 Further, the branches ${\cal C}_2$
and ${\cal C}_3$  are given by Theorem \ref{p1} and both of them
curve to the
right of $\lambda_1^+$ and  $\lambda_1^-$, respectively.
Neither ${\cal C}_2$ or ${\cal C}_3$  connect to ${{\cal C}_1}$, since
 ${\cal C}_1$ is isolated from the exterior of the
open set ${\cal O}$,
see Remark 5.2. Observe that  the elements of ${\cal C}_2$
(resp. ${\cal C}_3$) are outside ${\cal O}$ for $\lambda$ close to $\lpl$ (resp. $\lmin$).

Therefore the uniqueness statement of Proposition \ref{mininf} excludes the alternatives in Theorem \ref{p1} 3) (iii) and  4) (iii), since we already know that ${\cal C}_1$ contains solutions for arbitrary small $\lambda$. We already noted cases 3) (ii) and  4) (ii) are excluded by (F1). Finally, case 3) (iv) in Theorem \ref{p1} is excluded by the a priori bound in Proposition \ref{fi-}, so only case 3) (i) remains.
 \hfill
$\Box$

\medskip

\noindent
{\bf Proof of Theorem \ref{teo3}.}
We fix a small number $\varepsilon>0$ and for each $K>0$ consider a solution $u_K$ of $F[u_K]+ (\lmin-\varepsilon)u_K = K$ in $\Omega$, $u_K=0$ on $\partial \Omega$, $u_K<0$ in $\Omega$. We know such a function $u_K$ exists, by Theorem 1.9 in \cite{QS}. By (F0) we can fix $K_0$ such that $K_0>f(x,K_0)$ in $\Omega$, hence $\underline{u} = u_{K_0}$ is a sub-solution of \equ{first}, for all $\lambda\in (-\infty, \lmin-\varepsilon)$. The super-solution
to consider is $\overline u\equiv 0$, as given by hypothesis (F1).
Then  Proposition \ref{boundedbranch} yields the existence of a  branch
${\cal C}_1^{\varepsilon}$  such that
 ${\cal C}_1^{\varepsilon}(\lambda)\not =\emptyset$ for all $\lambda\in
(-\infty, \lambdaminus-\varepsilon)$.

Next we pass to the limit as $\varepsilon\to 0$, like in the proofs of
Proposition \ref{degree} and Theorem 1.5.1 in \cite{FQS}, and obtain
either a connected component of ${\cal S}$ which bifurcates from infinity
to the left of $\lmin$, or a {\it bounded} branch of solutions which
"survives" up to $\lmin$, and hence "continues" in some small right
neighborhood of $\lmin$, again like in the proof of Proposition \ref{degree}.
The first of these alternatives is
3. (i). In case the second alternative is realized  there is a connected set of solutions ${\cal C}_3$ bifurcating
from minus infinity towards the left of $\lambdaminus$, as predicted
  in Theorem \ref{p1}. We claim this branch contains only negative solutions. To prove this, we set $$A=\{ (\lambda,u)\in {\cal C}_3\:|\: \lambda\in(-\infty,\lmin),\; \max_\Omega u >0\}.$$
The set $A$ is clearly open in ${\cal C}_3$, and $A\not ={\cal C}_3$. Hence if $A$ is not empty, then $A$ is not closed in ${\cal C}_3$, by the connectedness of ${\cal C}_3$. This means there is a sequence  $(\lambda_n,u_n)\in A$ such that $\lambda_n \to \lambda$, $u_n\to u$, and the  limit function  $u$ satisfies $u\le 0$ in $\Omega$, $u$ vanishes somewhere in $\Omega$, and solves the equation $F[u] + (\lambda-c)u= f(x,u)-cu\ge 0$ in $\Omega$, for some large $c$. Hence by Hopf's lemma $u\equiv0$, a contradiction with (F1).

Therefore ${\cal C}_3$  cannot  connect with the
branch bifurcating from plus infinity at $\lambdaplus$. It is not connected to
${\cal C}_1$ either -- by the isolation property of ${\cal C}_1(\lambda)$, see Remark \ref{rem2}. Further, ${\cal C}_3$ cannot contain solutions for arbitrarily small~$\lambda$, since ${\cal C}_1$ does, and we know solutions are unique for sufficiently small $\lambda$. Hence ${\cal C}_3$ must eventually curve to the right, so extra solutions appear, proving 3. (ii) and 4.

Finally, a branch ${\cal C}_2$ bifurcating from plus infinity towards the
left of $\lambdaplus$ exists thanks to Theorem  \ref{p1}. This branch is
kept away from ${\cal C}_1$ and ${\cal C}_3$, as we already saw, and, again by the uniqueness of solutions for sufficiently small $\lambda$, ${\cal C}_2$ has to curve to the right. This completes the proof.\hfill
$\Box$

\medskip

The occurrence
of alternative 3. (ii) in Theorem  \ref{teo3} can be avoided if  $f$ satisfies some further hypotheses.

\begin{proposition}\label{propsuppl}
Under the hypotheses of Theorem \ref{teo3}, if  in addition we make one of the following assumptions
\begin{itemize}
\item[1.] $f(x,s)$ is concave in $s$ for $s<0$,
\item[2.] for each $a_0>0$ there exists $k_0>0$ such that
\begin{equation}\label{ferte}
\frac{f(x,-k\fipl)}{k}<f(x,-a\fipl), \qquad \mbox{for all }\; a\in(0,a_0), k>k_0,\end{equation}  \end{itemize}
 then  alternative 3. (ii) in Theorem \ref{teo3}
does not occur.
\end{proposition}

\noindent{\bf Remark.} Note the model example of a sub-linear nonlinearity which satisfies the hypotheses of Theorem \ref{teo3}
$$
f(x,s) = - s|s|^{\alpha-1} + h(x), \qquad \alpha\in(0,1), \;h\gneqq 0
$$
satisfies both hypotheses in the above proposition. \medskip

\noindent
{\bf Proof of Proposition \ref{propsuppl}.} We are going to prove the following stronger claim : under the hypotheses of the proposition,  there cannot exist sequences $\lambda_n$, $u_n$, $v_n$, such that $\lambda_{n}<\lambda_{n+1}$, $\lambda_n\to \lmin$, $u_n,v_n<0$ in $\Omega$, $\|u_n\|$ is bounded, $\|v_n\|\to \infty$ and $u_n$ and $v_n$ are solutions of \re{first} with $\lambda=\lambda_n$.

Assume this is false and 1. holds. Then (passing to subsequences if necessary) $u_n$ is convergent in $C^1(\overline{\Omega})$, and $v_n/\|v_n\|\to \fimin$ in $C^1(\overline{\Omega})$, so  there is $n_0$ such that for all $n\ge n_0$ we have
$v_{n}<u_{n+1}$ in $\Omega$. The negative function $u_{n+1}$ is clearly a strict subsolution of $F[u] +\lambda_n u = f(x,u)$, and, since the zero function is a strict supersolution of this equation, it has a negative solution which is above $u_{n+1}$. We define
$$\overline{v}_n=\inf\{ v\:|\: u_{n+1}<v<0, \; v\mbox{ is}\mbox{ a}\mbox{ supersolution}\mbox{ of } F[u] +\lambda_n u = f(x,u)\:\}.
$$
Then $\overline{v}_n$ is a solution of $F[u] +\lambda_n u = f(x,u)$ such that between $u_{n+1}$ and $\overline{v}_n$ no other solution of this problem exists. Indeed, $\overline{v}_n$ is a supersolution (as an infimum of supersolutions), so between $u_{n+1}$ and $\overline{v}_n$ there is a minimal solution, with which $\overline{v}_n$ has to coincide, by its definition. Note  Hopf's lemma trivially implies that for some $\varepsilon >0$ we have $v_n<u_{n+1}-\varepsilon \fipl <\overline{v}_n -2\varepsilon\fipl$.

Next, by the convexity of $F$ and the concavity of $f$ we easily check that the function $u_\alpha=\alpha v_n + (1-\alpha) \overline{v}_n$ is a supersolution of $F[u] +\lambda_n u = f(x,u)$, for each $\alpha\in[0,1]$. This gives a contradiction with the definition of $\overline{v}_n$, for $\alpha$ small enough but positive. \medskip

Assume now our claim is false and 2. holds. We again have $-C_0\fipl\le u_n \le -c_0\fipl<0$ and $v_n/\|v_n\|\to \fimin$, so the numbers
$$\varepsilon_n:=\sup \{ \varepsilon>0 \:|\: u_n\le \varepsilon v_n\mbox{ in }\Omega\,\} $$
clearly satisfy $\varepsilon_n>0$ and $\varepsilon_n\to 0$. Hypothesis \re{ferte} implies that for sufficiently large $n$ we have $\varepsilon_n f(x,v_n) <f(x,u_n)$, that is, $F[\varepsilon_n v_n] + \lambda_n \varepsilon_n v_n < F[u_n] + \lambda_nu_n$, and Hopf's lemma yields a contradiction with the definition of~$\varepsilon_n$.
\hfill  $\Box$

\setcounter{equation}{0}
\section{Discussion and examples }\label{examples}

The main point of this section is to provide some examples showing that
when $(F1)$ or $(F2)$ fail,  then the bifurcation diagram for \equ{first}
may look very
differently from what is described in Theorems \ref{teo2}-\ref{teo3}.
However, we begin with some general comments
on our hypotheses and their
use.

Hypothesis (F0) is classical sub-linearity for $f$, which guarantees
bifurcation from infinity and also  ensures the solutions
of \equ{first} tend to zero as $\lambda\to-\infty$.
Condition
(F1) guarantees the existence of  a strict super-solution of \equ{first} for
all $\lambda$, while (F2) is used in some comparison statements and to
prove uniqueness of solutions of \re{first} for sufficiently negative $\lambda$.

Further, conditions
 \hip-\him\ and
\hdp-\hdm\ are the Landesman-Lazer type hypotheses which give a priori
bounds when $\lambda$ stays on one side of the
eigenvalues, and thus provide a solution at resonance and determine on which side of each eigenvalue the
bifurcation from infinity takes place.
The strict inequalities in \hip-\hdm\ are important and cannot be relaxed in general - for instance the problem
$F[u] + \lpl u = -\sqrt{\max\{1-u,0\}}$ has no solutions (and hence Theorem \ref{teo1} fails), as Theorems 1.6 and 1.4 in \cite{QS} show, even though the nonlinearity satisfies
the hypotheses of Theorem \ref{teo1}, except for the strict inequality in \hip.
On the other hand, for $F=\Delta$ it is known that in the
case of equalities in \hip-\hdm\ one can give supplementary assumptions on $f$
and
the rate of convergence of $f$ to its limits $f_\pm,f^\pm$, so that
results like Theorem \ref{teo1} still hold, see for instance Remark 21 in
\cite{AG}. Extensions of these ideas to HJB operators are out of the scope of
this work and could be the basis of future research.


Now we discuss examples where $(F1)$ or $(F2)$ fail.

\medskip

\noindent{\bf Example 1}. Our first example shows that for all sufficiently small  $\delta>0$ we can construct   a nonlinearity $f$
which does not satisfy (F1) and for which the set $\cals(\lpl+\delta)$ is empty. This means that,
in the framework of Theorems \ref{teo1}-\ref{teo2}, the branch
bifurcating from infinity to the right of $\lpl$ "turns back" before it
reaches $\lpl+\delta.$
A similar situation can be described for the branch bifurcating from minus
infinity to the left of $\lmin$ that "turns right", before reaching
$\lmin-\delta$.
In particular there cannot be
a continuum of solutions along the gap between $\lpl$ and $\lmin$.

Consider the Dirichlet problem
\begin{equation}\label{refgt}
F[u]+\lambda u=t\fipl+h\;\mbox{ in }\; \Omega,\qquad u=0\;\mbox{ on }\; \partial\Omega,
\label{ex1}
\end{equation}
at resonance, that is, for $\lambda=\lpl$. When $t=\tplus(h)$
 equation \re{refgt} may or may not have a solution,
depending on $F$ and $h$. An example of such a situation was given in \cite{armstrong} and we recall it here.
Take
$F[u]
=\max\{\Delta u,2\Delta u\}$, and $h \in \co$  such that
$\int_\Omega h\varphi_1=0$ and $h$ changes sign on $\partial \Omega$.
Here
 $\lpl = \lambda_1$, $\lmin= 2\lambda_1$, $\fipl = -\fimin= \varphi_1$, where
$\lambda_1$ and  $\varphi_1$ are the first eigenvalue and eigenfunction
of the Laplacian.
Then (see Example 4.3 in \cite{armstrong})  under the above
hypotheses on $h$ we have $t^*_+=0$ and problem \re{ex1} has no solutions
if $\lambda=\lpl$ and $t=t^*_+$. By exactly the same reasoning it is possible to
show that problem \re{ex1} has no solutions if $\lambda=\lmin$ and $t=t^*_-$.
\begin{lemma}\label{lemm1}
If equation \equ{ex1} with $\lambda=\lpl$ and $t=t^*_+$ does not have a
solution then there exists
$\delta_0$ such that $\tl>\tplus$ provided $\lambda\in (\lpl,\lpl+\delta_0)$.
Similarly, if  \equ{ex1} with $\lambda=\lmin$ and $t=t^*_-$ does not have a
solution then there exists
$\delta_0$ such that  $\tl>\tminus$ whenever $\lambda\in (\lmin-\delta_0, \lmin)$.
\end{lemma}
Before proving the lemma, we use it to construct a nonlinearity with the desired
properties.
For $\lambda$ sufficiently close to
$\lpl$ we have $\tplus<\tl$ so that we can choose $\bar t\in (\tplus,\tl)$.
We then define
\begin{equation}\label{deff}
f(x,u) = \left\{\begin{array}{ccc} \bar t\fipl + h &\mbox{if}& u\ge - M \\
\left(\frac{\bar t-\tplus+\varepsilon}{M}(u+M)+\bar t\right)\fipl + h &\mbox{if}& -2M\le u\le - M\\
(\tplus-\varepsilon)\fipl + h &\mbox{if}& u\le - 2M,
\end{array}\right.\end{equation}
where $\varepsilon$ and $M$ are some positive constants.
We readily see that  $f$ satisfies (F0) and \hip, the hypotheses
of Theorem \ref{teo1}, but $\cals(\lambda)$ is empty.
Indeed, if $\overline{u}\in \cals(\lambda)$, then $\overline{u}$
is a  super-solution for \re{ex1}. On the other hand, by Theorem 1.9 in
\cite{QS}, the equation $F[u] + \lambda u = K\|\bar t\fipl + h\|_{\linf}$ with
$\lambda<\lmin$ has a solution $u_K$,
for each $K>0$. Moreover, for large $K$, $\underline{u}_K$ is a
sub-solution of \re{ex1} and  $\underline{u}_K < \overline{u}$.
Then by Perron's method \re{ex1} has a solution, a contradiction.

Similarly, for $\lambda<\lmin$ sufficiently close to $\lmin$ we choose $\bar t
\in (\tminus, \tl)$ and
define $f(x,u)$ being equal to $\bar t\fipl + h$ if $u\le M$ and to $(\tminus-
\varepsilon)\fipl + h$ if $u\ge 2M$. By the same reasoning we find that
$\cals(\lambda)$ is empty.

We summarize: with these choices of $\lambda$ and $f$ there is a region of non-existence in the gap between
$\lpl$ and $\lmin$.
In other words, the connected sets of
solutions  of \equ{first}
${\cal C}_2$ (resp. $\calc_3$), predicted in Theorem \ref{teo1},
 do not extend to the right (resp. to the left) of
$\lambda$.
The first graph at the end of this section is an illustration of this situation.

We observe that if we take $M$ sufficiently large then all solutions of \re{ex1}
and \re{first}, with $f$ as given in \re{deff}, coincide. In fact,
we can take $-M$ to be a lower bound for all solutions of the
inequality $F[u] + \lambda u\le c +h$, where $c$ is such that $f(x,u)\le c+h$
in $\Omega$. Such an $M$ exists by the one-sided ABP inequality given in
Theorem 1.7 in \cite{QS}. Now we see that \re{first} with this $f$ has a unique solution for $\lambda<\lpl$, as an application of Theorem 1.8 in \cite{QS}, and
then the branch of solutions bifurcating from plus infinity must turn left
and go towards infinity near the $\lambda$-axis, as drawn on the picture.

\noindent{\bf Proof of Lemma \ref{lemm1}}.
Given $\lambda\in(\lpl,\lmin)$, let
 $v^*_\lambda$ be a solution of
 \begin{equation}
F[u]+\lambda u=\tl\fipl+h\;\mbox{ in }\; \Omega,
\label{ex2}
\end{equation}
whose existence is guaranteed by the results in \cite{S1}. We notice that
 $\|v^*_\lambda\|$ is unbounded as $\lambda\searrow\lpl$,
as otherwise $v^*_\lambda$ a subsequence of
$v^*_\lambda$ would converge to a solution of \re{ex1} with
$\lambda=\lpl$ and $t=\tplus$, which is excluded by assumption.
That  $\tl\to \tplus$ as $\lambda\to \lpl$ was proved in
Proposition~\ref{frehg}.
Then,
by the simplicity of $\lpl$, we find that
$v^*_\lambda/\|v^*_\lambda\|_\infty\to\fipl$ as $\lambda\to\lambdaplus$,
 in particular, $v^*_\lambda$ becomes positive in $\Omega$,
for $\lambda$ larger than and close enough to  $\lambdaplus$.
Suppose for contradiction that $\tplus\ge t_\lambda^*$,  then $v^*_\lambda\ge0$
satisfies
$$
F(v^*_\lambda)+\lambdaplus v^*_\lambda\le F(v^*_\lambda)+\lambda v^*_\lambda
=t_\lambda^*\varphi_1^+ +h\le \tplus\varphi_1^+ +h,
$$
so $v^*_\lambda$ is a super-solution for \equ{ex1} with
$\lambda= \lambdaplus$ and $t=\tplus$. As we already showed above, \re{ex1}
has a sub-solution
below $v^*_\lambda$, providing a contradiction.

In the same way, we see that $v^*_\lambda/\|v^*_\lambda\|_\infty\to\fimin$ as $\lambda\nearrow\lambdaminus$ and then $v^*_\lambda$ becomes negative in $\Omega$,
for $\lambda<\lambdaminus$ and close enough to  $\lambdaminus$.
Then $\tminus\ge t_\lambda^*$ would imply that $v^*_\lambda\le0$ satisfies
$$
F(v^*_\lambda)+\lambdaminus v^*_\lambda\le F(v^*_\lambda)+\lambda v^*_\lambda
=t_\lambda^*\varphi_1^+ +h\le \tminus\varphi_1^+ +h,
$$
so $v^*_\lambda$ is a super-solution for \equ{ex1} with
$\lambda= \lambdaminus$ and $t=\tminus$. To construct a sub-solution we
consider  $v_\ep$ a solutions of $F[v_\ep] + \lmin v_\ep =
(t^*_-+\ep)\fipl + h$, with $\ep>0$. By our assumption,
$v_ \ep/\|v_\ep\|_\infty\to\fimin$ as $\ep\to0$ (see  also Theorem 1.4 in \cite{FQS}).
Hence
there exists $\ep=\ep(\lambda)$ such that $v_\ep < v^*_\lambda$
and $F[v_\ep] + \lmin v_\ep \ge\tminus\varphi_1^+ + h$ and then Perron's method gives a contradiction again. \hfill $\Box$

\bremark
The Claim gives an idea of the behavior of $t_\lambda^*$, with respect to $\lambda$, near the extremes of
the interval $[\lpl,\lmin]$. However we do not have any idea about the global
behavior of $t_\lambda^*$, actually we do not even know how   $t_+^*$ and $t_-^*$ compare.
\eremark

For completeness  we  give a direct proof of the fact that in the above
examples  condition (F1) is not satisfied by  nonlinearities like in
\re{deff}. In this direction we have the following lemma, which is
 of independent interest.

\begin{lemma}\label{lemm2} For any $h\in \lp$, $p>N$, which is not a multiple of $\fipl$,
\begin{itemize}
\item[(a)]\ if $h\ge0$ and $h\not \equiv 0$ then $\tplus(h)<0$ and $\tminus(h)<0$;
    \item[(b)]\ if $h\le0$ and $h\not \equiv 0$ then $\tplus(h)>0$ and $\tminus(h)>0$;
\item[(c)]\ the functions $\tplus(h)\fipl + h$ and $\tminus(h)\fipl + h$ change sign in $\Omega$. \end{itemize}
    \end{lemma}

\noindent{\bf Proof.} {\it (a)} If  $\tplus(h)\ge 0$ then, as $h\ge 0$, by Theorem 1.9 in
\cite{QS} the problem $F[u] + \lpl u = \tplus(h)\fipl + h
$ has a solution. Then by Theorem 1.2 in \cite{FQS} $u+k\fipl$ is a solution of the
same problem, for all $k>0$. Since $u+k\fipl$ is positive for sufficiently large $k$, by Theorem 1.2 in \cite{QS} we get that $u$ is a multiple of $\fipl$, a contradiction, since $h\not =0$.

If $\tminus(h)\ge 0$, by Theorem 1.5 in \cite{FQS} either there exist
sequences $\ep_n\to 0$ and $u_n$ of solutions of the problem
$F[u_n] + \lmin u_n = (\tminus(h)+\ep_n)\fipl + h
$  such that $u_n$ is unbounded and $u_n$ is negative for large $n$,
or $F[u+k\fimin] + \lmin (u+k\fimin) = \tminus(h)\fipl + h
$ for some $u$ and all $k>0$. In both cases we get a negative solution of $F[u]+\lmin u \gneqq 0$, which by Theorem 1.4 in \cite{QS}  is then  a multiple of $\fimin$, a contradiction.

\noindent {\it (b)} If $\tplus(h)\le 0$ then $F[u] + \lpl u = \tplus(h)\fipl
+ h
$ has no solution by
Theorems 1.6 and 1.4 in \cite{QS}, since $\tplus(h)\fipl
+ h
 \lneqq 0$.
If $\tminus(h)\le 0$ we again have $\tminus(h)\fipl + h
\lneqq 0
$, then    $F[u] + \lmin u =
\tminus(h)\fipl + h
$ has no solutions by the anti-maximum principle, see for instance
Proposition 4.1 in \cite{FQS}.
Hence by Theorems 1.2 and 1.4 in \cite{FQS} there exist sequences
$\ep_n\to 0$, $u_n^+$ and $u_n^-$ of solutions of $F[u_n^\pm] + \lambda^\pm_1 u_n^\pm = (t^*_\pm(h)+\ep_n)\fipl + h  $
such that $u_n^\pm/ \|u_n^\pm\|_\infty \to \varphi^\pm$.
Fix $w$ to be the solution of the Dirichlet problem $F(w)-\gamma w = -h
$ in $\Omega$.
This problem is uniquely solvable, with $w<0$ in $\Omega$, since by (H1) the
operator $F-\gamma$ is decreasing in $u$ (see for instance \cite{CCKS} and
\cite{QS}).  Then by the maximum principle and Hopf's lemma
$\ep_n\fipl + (\lpl+\gamma) w<0$ in $\Omega$, if $n$ is sufficiently large.
Hence $u_n^++ w$ is positive and $F[u_n^++ w] + \lpl(u_n^++ w)<0$ in $\Omega$,
 which is a contradiction with Theorem 1.4 in \cite{QS}. Similarly,
$u_n^-+ w$ is negative and satisfies $F[u_n^-+ w] + \lmin(u_n^-+ w)<0$ in $\Omega$, which is a contradiction with Theorem 1.2 in \cite{QS}.

\noindent {\it (c)} This is an immediate consequence of (a) and (b). Indeed,
 if (c) is false we just replace $h$ by $t^*_\pm(h)\fipl + h$ in (a) or (b). \hfill $\Box$

\bremark
The statements on $\tplus$ in the preceding lemma also follow from Theorem 1.1 and formula (1.12) in \cite{armstrong}.
\eremark

The following example illustrate the role of hypothesis (F2), which allows the
 use of the  method of sub- and super-solutions, and prevents the branches
which bifurcate from
infinity to survive for arbitrarily negative $\lambda$.

\medskip

\noindent {\bf Example 2.}
Consider the function  $\omega(u)=\displaystyle \frac{u}{\sqrt{|u|}}$,
$\omega(0)=0$ and  the problem
\begin{equation}\label{simpi}
    \Delta u + \lambda u = -\omega(u), \;\mbox{ in }\Omega.
\end{equation}
This problem is variational and its associated functional is
$$
J(u) = \int_\Omega \left(|\nabla u|^2 -\lambda u^2 - |u|^{3/2}\right) \, dx,
$$
which is even, bounded below, takes negative values and attains its minimum on
$H_0^1(\Omega)$, for each $\lambda<\lambda_1$. The same is valid for
$J_+(u) = J(u^+)$ and $J_-(u) = J(u^-)$, whose minima are then a positive and a negative solution of~\re{simpi}.

In the context of nonlinear HJB operators we may consider
\begin{equation}\label{freza}
\max\{\Delta u, 2\Delta u\} + \lambda u= -\omega(u), \;\mbox{ in }\Omega\qquad u=0 \;\mbox{ on }\partial \Omega.
\end{equation}
For this problem we have bifurcation from plus infinity to the left of $\lambda_1$
and from minus infinity to the left of $2\lambda_1$.
These branches cannot reach the trivial solution set $\R\times \{0\}$, since
bifurcation of positive or negative solutions from the trivial solution does not occur for \equ{simpi}. Exactly as in the proof of Theorem \ref{teo3} (see the definition of the set $A$ in the previous section) we can show that they contain only positive or negative solutions.
Actually these branches are curves which can never turn, since positive and negative
solutions of  \re{freza} are unique -- this can be proved in the same way as Proposition~\ref{uniqi} below.


\noindent{\bf Example 3.} Finally, let us look at an example of a sub-linear nonlinearity $f$ which satisfies (F2) but $f(x,0)\equiv 0$. For any HJB operator $F$ satisfying our hypotheses consider
\begin{equation}\label{simps}
    F[u] + \lambda u = \tilde f(u):= \left\{ \begin{array}{ccr}
                     -u &\mbox{if}& |u|\le 1\\
                     -\omega(u)&\mbox{if}& |u|\ge 1\end{array}\right.
\end{equation}
In this situation we have positive (resp. negative) bifurcation from zero at $\lambda= \lpl-1$ (resp. $\lambda= \lmin-1$), more precisely, $(\lpl-1,k\fipl)$ and $(\lmin-1,k\fimin)$ are solutions for $k\in [0,1]$ (for more general results on bifurcation from zero see \cite{BEQ}). Further, note that there are only positive (resp. negative) solutions on these branches, as well on the  branches which bifurcate from plus (resp. minus) infinity, given by Theorem \ref{teo3}. This is a simple consequence of the strong maximum principle and the fact that the right hand side of \re{simps} is positive (resp. negative) if $u$ is negative (resp. positive), so if $u\le(\ge) 0$ and $u$ vanishes at one point then $u$ is identically zero. The bifurcation branches connect, as shown by the following uniqueness result.

\begin{prop} \label{uniqi} If $u$ and $v$ are two solutions of \re{simps} having the same sign and $\|u\|>1$ or $\|v\|>1$ then $u\equiv v$. If $\|u\|\le1$ and $\|v\|\le1$ then (by the simplicity of the eigenvalues) $\lambda = \lambda_1^\pm - 1$ and $u = v + k \varphi_1^\pm$ for some $k\in [0,1]$.
\end{prop}

\noindent{\it Proof.} Say $u>0$, $v>0$, $\|v\|>1$. Set
$$
\tau :=\sup \{\mu>0\;|\; u\ge \mu v\;\mbox{ in }\Omega\}.
$$
By Hopf's lemma $\tau>0$ and we have $u\ge \tau v$.

First, suppose $\tau<1$. By the definition of $\tilde {f}$ in \re{simps} and $\|v\|>1$ we easily see that
$$
\tilde {f} (u) \le \tilde {f} (\tau v) \lneqq \tau \tilde {f}(v)\;\mbox{ in }\Omega.
$$
Hence \re{simps} and the hypotheses on $F$ imply
$$
\mm(D^2(u-\tau v)) - \gamma |D(u-\tau v)| -(\gamma +  \lambda)(u-\tau v)\lneqq 0 $$
and $ u-\tau v\ge 0$  in $\Omega$, so Hopf's lemma implies $u\ge (\tau + \varepsilon)v$ for some $\varepsilon>0$, a contradiction with the definition of $\tau$.

Second, if $\tau \ge 1$ we repeat the above argument with $u$ and $v$ interchanged. This leaves $u\ge v $ and $ v\ge u$ as the only case not excluded.\hfill $\Box$

\medskip

The following picture summarizes  the above examples.
\begin{figure}[h]
\includegraphics[width=14cm, height=4cm]{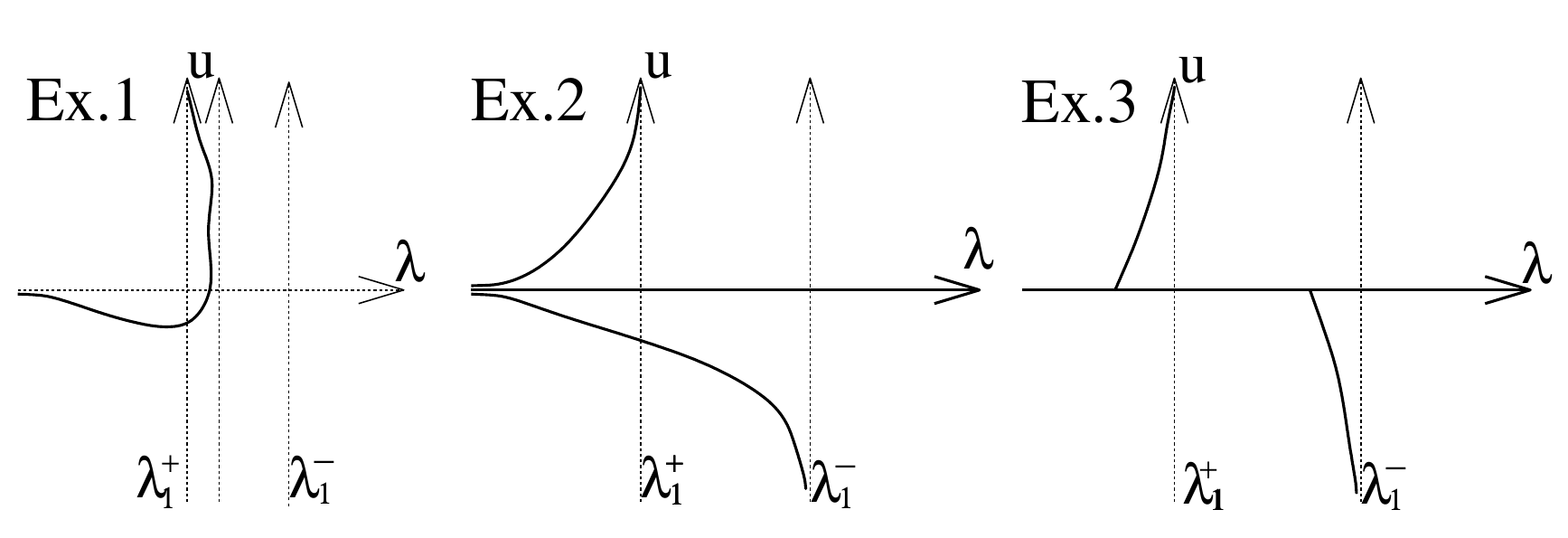}\label{fig2}
\end{figure}

\bigskip

\noindent
{\bf Acknowledgements:}
P.F. was  partially supported by Fondecyt Grant \# 1070314,
FONDAP and BASAL-CMM projects
 and  Ecos-Conicyt project C05E09.
A. Q. was partially supported by Fondecyt Grant \# 1070264 and USM Grant \#   12.09.17
and Programa Basal, CMM. U. de Chile. .

\begin{flushleft}
Patricio FELMER\\ Departamento de
Ingenier\'{\i}a  Matem\'atica,
Universidad de Chile, Casilla 170 Correo 3,
  Santiago, Chile.\\
  e-mail : pfelmer@dim.uchile.cl
  \bigskip

Alexander QUAAS\\
Departamento de  Matem\'atica,
Universidad T\'ecnica Federico Santa Mar\'{i}a,
Casilla: V-110, Avda. Espa\~na 1680,
Valpara\'{\i}so, Chile.\\
e-mail : alexander.quaas@usm.cl
\bigskip

Boyan SIRAKOV (corresponding author) \\
UFR SEGMI, Universit\'e de Paris 10,
92001 Nanterre Cedex, France, and
CAMS, EHESS,
54 bd. Raspail,
75006 Paris, France\\
e-mail : sirakov@ehess.fr

\end{flushleft}
\end{document}